\documentclass[pdflatex,sn-mathphys-num]{sn-jnl}

\usepackage[utf8]{inputenc}
\usepackage[T1]{fontenc}
\usepackage{lmodern}
\usepackage{booktabs}
\usepackage{amsmath,amssymb,amsfonts}
\usepackage{nicefrac}
\usepackage{microtype}
\usepackage{graphicx}
\usepackage{tabularx}
\let\orcidlogo\relax
\usepackage{orcidlink}

\usepackage{algorithm}
\usepackage{algorithmicx}
\usepackage{algpseudocode}

\usepackage{amsthm}
\theoremstyle{thmstyleone}

\newtheorem{proposition}{Proposition}

\usepackage[title]{appendix}
\usepackage{placeins}

\usepackage{subcaption}

\DeclareCaptionLabelSeparator{snfig}{\hskip.7em}
\DeclareCaptionLabelSeparator{sntab}{\hskip2mm}
\hypersetup{bookmarksdepth=3}

\DeclareMathOperator*{\argmin}{argmin}

\newcommand{\panelfigwidth}{178.56bp}
\newcommand{\singlefigwidth}{238.08bp}
\newcommand{\hspacegapsize}{0.015}

\begin{document}

\title[Incremental SVD for Large-Scale Dynamic Matrices]{Incremental SVD for Large-Scale Dynamic Matrices: Accuracy, Subspace Stability, Refresh Strategies, and 
Financial Factor-Based Risk Models}

\author*[1]{\fnm{Stilyan} \sur{Staykov}\,\orcidlink{0009-0006-2193-033X}}\email{s.staykov@student.tue.nl}

\affil*[1]{\orgdiv{Department of Mathematics and Computer Science}, \orgname{Eindhoven University of Technology}, \orgaddress{\city{Eindhoven}, \country{The Netherlands}}}

\abstract{Return panels, covariances, and large feature matrices evolve one observation or entry at a time, yet downstream models require an up-to-date low-rank factorization $A_t \approx U_t \Sigma_t V_t^\top$ on every tick---a regime where a full singular value decomposition (SVD) is prohibitive and existing alternatives sacrifice singular vectors, singular values, or long-horizon stability.
We present a metric-driven study of Brand-style \emph{incremental SVD}, built around a unified engine that handles row appends, column appends, rank-1 entry updates, and metrics tracking.
For rank-1 entry updates, we derive an explicit projection-based rule $U'\Sigma'(V')^\top = P_U(\widehat{A} + \delta\,e_ie_j^\top)P_V$ that keeps rank fixed while discarding only the out-of-subspace remainder quantifiably, turning Brand's rank-suppression heuristic into an operational scheme.
Each update discards a part whose norm is available at negligible cost; the error at the last refresh plus the sum of these norms since then is an a posteriori error bound $\beta_t$.
We then treat \emph{refresh scheduling} as a first-class design axis, comparing periodic refresh with a bound-triggered rule that refreshes once $\beta_t/\|A_t\|_F$ exceeds a tolerance.
The engine tracks error ratios, principal angles, explained variance, and per-update runtime on long synthetic streams and a factor model of exchange-traded funds for covariance and portfolio-risk estimation. With a sensible rank and refresh cadence, incremental SVD matches full-SVD accuracy within a few percent at a fraction of the cost. The bound-triggered rule keeps the relative error below any tolerance above its rank-$k$ optimum, with accuracy close to periodic refresh at a matched number of full SVDs.}



\maketitle

\section{Introduction and Background}
\label{sec:intro}
\label{sec:background}

\subsection{Problem}

Many data matrices of interest are not static objects but evolving streams. In recommender systems, user--item rating matrices change continually as new users arrive, new items are added, and existing users update their feedback. In quantitative finance, time$\times$asset return panels grow one row at a time as markets generate new observations and the cross-sectional universe of tradable instruments shifts.  In both settings, downstream tasks---such as rating prediction or risk modeling---are often built on top of a low-rank representation of the data, for example via a truncated singular value decomposition (SVD) or a closely related matrix factorization.

We briefly review SVD and low-rank approximation, setting up notation used throughout the paper.

\label{sec:background-svd}
For a matrix $A \in \mathbb{R}^{m \times n}$ the singular value decomposition (SVD) writes
$
A = U \Sigma V^\top,
$
where $U \in \mathbb{R}^{m \times m}$ and $V \in \mathbb{R}^{n \times n}$ are orthogonal, and
$\Sigma = \mathrm{diag}(\sigma_1,\dots,\sigma_r) \in \mathbb{R}^{m \times n}$ contains the singular values
$\sigma_1 \geq \sigma_2 \geq \dots \geq \sigma_r > 0$ with $r = \mathrm{rank}(A)$.
The truncated SVD of rank $k$ keeps only the leading $k$ singular values and vectors,
$
A_k = U_k \Sigma_k V_k^\top,
$
with $U_k \in \mathbb{R}^{m \times k}$, $V_k \in \mathbb{R}^{n \times k}$, and $\Sigma_k \in \mathbb{R}^{k \times k}$.
Throughout the paper, $\|\cdot\|$ denotes the 2-norm and $\|\cdot\|_F$ the Frobenius norm.
By the Eckart–Young–Mirsky theorem, $A_k$ is the best rank-$k$ approximation to $A$ in any unitarily invariant norm, in particular in the Frobenius norm,
\begin{equation}
A_k \in \argmin_{\mathrm{rank}(B)\le k} \|A - B\|_F.
\label{eq:eckart-young-frob}
\end{equation}
The squared Frobenius norm decomposes as
$\|A\|_F^2 = \sum_{i=1}^r \sigma_i^2$ and the residual error of $A_k$ equals
$\|A - A_k\|_F^2 = \sum_{i=k+1}^r \sigma_i^2$.
This naturally leads to the explained variance ratio (EVR), defined by
$
\mathrm{EVR}(k) = \frac{\|A_k\|_F^2}{\|A\|_F^2}
= \frac{\sum_{i=1}^k \sigma_i^2}{\sum_{i=1}^r \sigma_i^2}
= 1 - \frac{\|A - A_k\|_F^2}{\|A\|_F^2},
$
\label{eq:EVRk}
which measures how much of the total variability is captured by the leading $k$ singular values.
In static settings one typically chooses $k$ such that $\mathrm{EVR}(k)$ exceeds a target
(e.g.\ $90$–$95\%$), or by inspecting the decay of $\sigma_i$.

In our work, SVD and its truncation play two roles.
First, the truncated SVD $U_t \Sigma_t V_t^\top$ of a dynamic matrix $A_t$ is the model we maintain incrementally.
Second, the full SVD of $A_t$ provides an \emph{oracle baseline}:
its truncation $A_t^{\text{opt}}$ defines the optimal rank-$k$ error against which we compare the incremental approximation.
Section~\ref{subsec:metrics} formalizes these metrics.

\subsection{Motivation}

Recomputing a full SVD after every update is typically infeasible in systems that operate under latency and throughput constraints.  This motivates \emph{incremental} or \emph{online} SVD methods~\cite{brand2006fast}, which update a rank-$k$ approximation $A_t \approx U_t \Sigma_t V_t^\top$ as the matrix $A_t$ evolves, for example via row/column appends or rank-1 entry updates.

\label{sec:background-finance}
A concrete instance arises in equity and multi-asset finance. One constructs log returns $r_{t,i} = \log(P_{t,i} \, / \,P_{t-1,i})$ from asset prices $P_{t,i}$, yielding a matrix $R_T \in \mathbb{R}^{T \times N}$ of $T$ days across $N$ assets.
PCA---equivalent to SVD on mean-centered returns---is a standard tool for extracting systematic risk factors: the leading right singular vectors of $R_T$ give cross-sectional factors, while the associated singular values determine their variance contribution.

Low-rank factor models underlie many portfolio construction and risk management workflows: a rank-$k$ factorization $R_T \approx F_T B^\top$ (with $F_T = U_k\Sigma_k$ and $B = V_k$) yields a low-rank-plus-diagonal covariance that can be propagated through portfolio weights to estimate volatility and Value-at-Risk.

Modern trading and risk systems operate online: returns arrive continuously and covariance and factor estimates must be updated under tight latency constraints.

\subsection{State of the Art}
\label{sec:background-incremental-svd}

When new rows, columns, or localized updates arrive over time, recomputing a full SVD is often too expensive, motivating incremental, online, and streaming variants.
Brand~\cite{brand2006fast,brand2002incremental} gives efficient update formulas for augmenting a (possibly truncated) SVD by a row, column, or low-rank perturbation: the update is expressed in the current singular-vector basis, producing a small ``core'' matrix whose SVD captures the change to the spectrum.

Complementary update backbones have since been developed: projection-based partial-SVD updates via Rayleigh--Ritz procedures in the LSI setting~\cite{vecharynski2014fastlsi}; fast truncated-SVD updates for large sparse matrices in representation learning~\cite{deng2024fasttsvd}; subspace-tracking with explicit attention to canonical angles and long-horizon stability~\cite{moonen1992svdtrack}; and recent ``SVD-type'' streaming updates that revisit accuracy--throughput trade-offs via related factorizations~\cite{brust2025svdtype}. More recent projection-based schemes refine this idea by constructing approximation subspaces that match the leading singular vectors of the augmented matrix, offering improved accuracy on latent semantic indexing and recommender-system datasets~\cite{kalantzis2021projection}.

In finance specifically, existing online approaches---rolling-window PCA, dynamic conditional correlation (DCC) models, Kalman filter-based state-space methods---mostly recompute PCA per window or model covariance dynamics parametrically; few directly study exact incremental SVD updates on the growing returns matrix.

\subsection{Limitations of the State of the Art}

Despite this rich literature, several practical challenges remain under-explored:
(i)~how approximation error and subspace drift accumulate over long streams;
(ii)~how to choose and adapt the working rank $k$ in a streaming regime; and
(iii)~how to design policies for \emph{refreshing} the model, i.e., when to pay the cost of recomputing a full SVD in order to re-align with the current ``true'' low-rank structure.  These questions are particularly acute in applications where the latent factors themselves have a semantic interpretation, such as user/item embeddings in recommenders or risk factors and covariance structures in finance.

More broadly, most works either focus on theoretical guarantees under idealized noise, or present small empirical examples without tracking how error, drift, and runtime evolve over long streams; the design of \emph{refresh policies}---when to recompute a full SVD to rebaseline---is usually left implicit.

Classical SVD and PCA provide the optimal low-rank approximation in static settings, and finance relies heavily on low-rank factor models.
Incremental and online SVD methods exist, but practical aspects such as error accumulation, subspace drift, rank choice, and refresh policies are rarely analyzed systematically in the context of online risk models.
This motivates our metric-driven study of incremental SVD for dynamic matrices in Sections~\ref{sec:incremental-svd}--\ref{sec:finance}.

\subsection{Contributions}

Our work takes a complementary, mainly empirical view: given a concrete Brand-style implementation, we bound the error that accumulates between refreshes and study how rank choice and refresh strategies affect accuracy, stability, and runtime in synthetic and applied settings. This paper takes a deliberately practical view on incremental SVD for dynamic matrices, with an emphasis on accuracy, stability, and refresh strategies in settings motivated by recommender systems and financial factor models. Our main contributions are:

\begin{enumerate}
  \item \textbf{Unified incremental SVD engine for dynamic matrices.}
    We implement and analyze a single incremental SVD framework that supports
    three fundamental update types for a time-varying matrix $A_t$:
    row appends (new users or new time steps), column appends (new items or
    new assets), and rank-1 entry updates (individual ratings or returns).
    The updates are \emph{exact up to rank-$k$ truncation}, and can be
    combined with periodic ``refresh'' steps that recompute a truncated SVD
    of the current matrix.  This yields a flexible engine that mirrors the
    operations encountered in realistic recommender and finance pipelines.

  \item \textbf{A posteriori error control between refreshes.}
    We prove that the part of the matrix discarded by each update has a
    norm that is available in closed form or as a by-product of the
    update, and that the error at the last refresh plus the sum of these
    norms since then is a bound $\beta_t$ on the incremental reconstruction error
    (Proposition~\ref{prop:error-control}).  The engine maintains
    $\beta_t$ at negligible cost per step and without a full SVD; on
    segments with a single update type the error itself follows from
    $\|A_t\|_F$ and the maintained singular values.  The bound holds at
    every step of our experiments, but it is loose on long segments
    without refresh.

  \item \textbf{Metric-driven framework for accuracy and subspace stability.}
    We propose a unified set of metrics to evaluate incremental SVD:
    Frobenius reconstruction error versus the optimal rank-$k$ error
    (from a full SVD), error gap and error ratio, principal angles between
    the maintained subspace and (i)~its initial basis and (ii)~the current
    optimal basis, explained variance ratio, and runtime statistics for both
    incremental updates and full SVD recomputations.  This framework makes
    it possible to quantify not only how well the approximation tracks the
    best rank-$k$ SVD, but also how the latent subspace drifts over time and
    what computational trade-offs different policies incur.

  \item \textbf{Refresh and rank-selection policy design space.}
    Building on this engine and metric set, we formalize a small
    but expressive design space of policies for controlling drift and model
    complexity.  These include simple periodic refresh, error-based refresh
    (triggered by an increase in error ratio), angle-based refresh
    (triggered by principal-angle thresholds), and EVR- or residual-based
    adaptive rank rules.  Our experiments compare periodic refresh with a
    bound-triggered rule that refreshes once $\beta_t/\|A_t\|_F$ exceeds a
    tolerance $\tau$, on a synthetic rank-1 stream and on the ETF stream.
    As long as $\tau$ exceeds the relative error of the best rank-$k$
    approximation, the rule keeps the reconstruction error below
    $\tau\|A_t\|_F$ after every step; at a matched number of full SVDs its
    accuracy is close to that of periodic refresh.  The recorded error ratios, principal
    angles, and explained variance provide the signals for the angle-based
    and adaptive-rank rules.

  \item \textbf{Application case study in
  finance.}
    Finally, we instantiate the framework on a representative application.
    On a panel of multi-asset exchange-traded funds, we interpret the SVD as an
    online factor model, and evaluate how well incremental SVD tracks the
    leading eigenmodes of the return covariance and the associated portfolio
    risk for test portfolios.  This case study illustrates how the same
    incremental machinery can be used for
    factor-based risk modeling in a way that is compatible with standard
    datasets and data providers.
\end{enumerate}

Our finance experiments build a returns matrix for a universe of liquid ETFs, maintain an incremental SVD as new days arrive, and compare the resulting low-rank covariance and portfolio risk estimates to full-SVD recomputation; the same pipeline extends to institutional sources such as WRDS (e.g.\ CRSP or TAQ).

\subsection{Paper Overview}

Section~\ref{sec:incremental-svd}
formalizes the dynamic matrix setting, describes the incremental update rules
for rows, columns, and rank-1 entry updates, and defines the error and
subspace metrics we use throughout. Section~\ref{sec:policies} introduces
refresh and adaptive rank policies. Section~\ref{sec:synthetic} presents
synthetic streaming experiments that probe accuracy, subspace drift, and
runtime scaling.  Section~\ref{sec:finance} applies the framework to
a multi-asset ETF return panel. Finally, Sections~\ref{sec:discussion} 
and~\ref{sec:conclusion} conclude the paper and discuss limitations and directions for future work.

\section{Incremental SVD Framework}
\label{sec:incremental-svd}

Now we formalize the dynamic matrix setting, derive the incremental update rules implemented in our \texttt{IncrementalSVD} engine, and define the error and subspace-drift metrics that will be used throughout the experiments. The goal is to make explicit what the algorithm does for each type of update (rows, columns, and single entries), how it interacts with periodic full-SVD ``refreshes'', and what quantities we monitor to quantify accuracy and stability.

Our framework builds on the deterministic Brand-style updates~\cite{brand2006fast,brand2002incremental}, which maintain an \emph{exact} rank-$k$ factorization after each append or rank-1 perturbation.
We refer to~\cite{balzano2018streaming} for a unified survey of streaming PCA and subspace-tracking algorithms, organized by algebraic and geometric perspectives and reviewing both asymptotic and non-asymptotic guarantees, which complements the algorithm-by-algorithm comparison below.
Alternative online paradigms include \emph{randomized SVD}~\cite{halko2011randomized}, which is efficient for batch recomputation but does not natively support incremental updates; \emph{matrix sketching} methods such as Frequent Directions~\cite{liberty2013simple,ghashami2016frequent}, which absorb streaming rows but do not maintain explicit singular vectors or support entry-level updates; and \emph{stochastic subspace trackers} like Oja's rule~\cite{oja1982simplified} and GROUSE~\cite{balzano2010grouse}, which scale to high dimensions but offer only asymptotic guarantees and do not produce singular values.
The Brand-style approach produces a full $(U, \Sigma, V)$ triple at every step, which is essential for our metric-driven evaluation; we revisit these trade-offs in Section~\ref{sec:discussion}.

\subsection{Problem Setting and Notation}

We consider a real-valued matrix $A_t \in \mathbb{R}^{m_t \times n_t}$ that evolves over time for $t=0,1,2,\dots$, where rows and columns may be appended and individual entries may be modified. At any time $t$ we maintain a rank-$k$ truncated singular value decomposition $A_t \approx U_t \Sigma_t V_t^\top$, where $U_t \in \mathbb{R}^{m_t \times k}$ and $V_t \in \mathbb{R}^{n_t \times k}$ have orthonormal columns, and $\Sigma_t = \mathrm{diag}(\sigma_{t,1},\dots,\sigma_{t,k}) \in \mathbb{R}^{k \times k}$ contains non-increasing singular values. Depending on the policy chosen at implementation time, the dimensions $(m_t, n_t)$ may be kept static---so that newly arriving information is absorbed into a fixed-size representation---or grown dynamically according to a criterion that decides when an actual size increase is warranted. We discuss both regimes; in our experiments the dimensions either stay fixed or grow with every append.

We explicitly model three types of updates:
\begin{enumerate}
    \item \textbf{Row append (new ``user'' / new time point).} A new row $x^\top \in \mathbb{R}^{1 \times n_t}$ arrives and the matrix becomes
    \begin{equation}
    A_{t+1} =
    \begin{bmatrix}
        A_t \\ x^\top
    \end{bmatrix}
    \in \mathbb{R}^{(m_t+1) \times n_t}.
    \label{eq:rowAppendForm}
    \end{equation}

    \item \textbf{Column append (new ``item'' / new asset).} A new column $y \in \mathbb{R}^{m_t}$ arrives and the matrix becomes
    \[
    A_{t+1} =
    \begin{bmatrix}
        A_t & y
    \end{bmatrix}
    \in \mathbb{R}^{m_t \times (n_t+1)}.
    \]

    \item \textbf{Single entry modification (rank-1 update).} A single entry at position $(i,j)$ is modified by $\delta$,
    \begin{equation}
    A_{t+1} = A_t + \delta\, e_i e_j^\top,
    \label{eq:rank1UpdateForm}
    \end{equation}
    where $e_i$ and $e_j$ are standard basis vectors.
\end{enumerate}

The core of our framework is an incremental algorithm described by~\cite{brand2002incremental,brand2006fast} that updates $(U_t,\Sigma_t,V_t)$ after each such operation, while keeping the working rank bounded by a user-specified $k$. In addition, we allow occasional \emph{full refreshes}, where a truncated SVD of the current matrix is recomputed from scratch. These refreshes serve both as a numerical ``reset'' and as a reference for measuring drift.

Algorithm~\ref{alg:incrementalsvd} gives the engine control flow; detailed pseudocode for each update primitive
appears in Appendix~\ref{app:updates} (Algorithms~\ref{alg:rowappend}--\ref{alg:rank1}). Note, in Sections~\ref{subsec:row-col-updates} and~\ref{subsec:rank1-entry-updates} we discuss the mathematical logic behind the algorithm.

\begin{algorithm}[t]
\caption{\texttt{IncrementalSVD} engine (high level). Detailed update routines are in Appendix~\ref{app:updates}}
\label{alg:incrementalsvd}
{\footnotesize
\begin{algorithmic}[1]
\Require Initial matrix $A_0$, target rank $k$, tolerance \texttt{tol}, refresh rule $\pi$
\State $(U,\Sigma,V^\top) \gets \texttt{TruncatedSVD}(A_0,k)$
\State $A \gets A_0$
\For{$t=1,2,\ldots$}
  \State Receive update event $\mathcal{E}_t$
  \If{$\mathcal{E}_t=\texttt{RowAppend}(x^\top)$}
    \State Append $x^\top$ to $A$ as the bottom row
    \State $(U,\Sigma,V^\top) \gets \texttt{RowAppendUpdate}(U,\Sigma,V^\top,x^\top;\texttt{tol})$
  \ElsIf{$\mathcal{E}_t=\texttt{ColAppend}(y)$}
    \State Append $y$ to $A$ as the right-most column
    \State $(U,\Sigma,V^\top) \gets \texttt{ColAppendUpdateViaTranspose}(U,\Sigma,V^\top,y;\texttt{tol})$
  \ElsIf{$\mathcal{E}_t=\texttt{RankOne}(i,j,\delta)$}
    \State $A_{i,j} \gets A_{i,j} + \delta$
    \State $(U,\Sigma,V^\top) \gets \texttt{RankOneUpdate}(U,\Sigma,V^\top,i,j,\delta)$
  \EndIf
  \State \textbf{Truncate:} keep the top-$k$ singular values and corresponding vectors
  \If{$\pi(t,\text{metrics})=\textsc{true}$}
    \State $(U,\Sigma,V^\top) \gets \texttt{TruncatedSVD}(A,k)$ \Comment{refresh from scratch}
  \EndIf
\EndFor
\end{algorithmic}
}
\end{algorithm}

In the synthetic experiments of Section~\ref{sec:synthetic} we work with dense matrices and measure error in the full Frobenius norm. 

\subsection{Incremental Row and Column Updates}
\label{subsec:row-col-updates}

We first describe how we update the SVD when appending a new row. Our implementation follows the small-core update pattern of Brand-style algorithms~\cite{brand2002incremental}: the effect of the new row is compressed into a $(k+1) \times (k+1)$ matrix, on which we compute an exact SVD. Column additions are handled with the transpose trick. 

Related projection-based truncated-SVD updating schemes, particularly in latent semantic indexing, provide a complementary Rayleigh--Ritz viewpoint and strategies for reducing the dimension of the update subspace~\cite{vecharynski2014fastlsi}.

\subsubsection{Row Append}

Assume we have a rank-$r$ truncated SVD
$
A_t \approx U \Sigma V^\top,\quad
U \in \mathbb{R}^{m \times r},\ V \in \mathbb{R}^{n \times r}.
$
We append a new row $x^\top \in \mathbb{R}^{1 \times n}$ to form the matrix
\[
A' =
\begin{bmatrix}
A_t \\ x^\top
\end{bmatrix}
\in \mathbb{R}^{(m+1) \times n}.
\]
The algorithm proceeds in three steps:

\paragraph{Step 1: projection onto current right-singular subspace.}
We decompose the new row into a part that lies in the span of $V$ and an orthogonal residual. Let
\begin{equation}
p = x^\top V, \qquad z = x^\top - p \, V^\top, \qquad \rho = \|z\|.
\label{eq:row-residual}
\end{equation}
If $\rho > \mathrm{tol}$, where \texttt{tol} is a small numerical threshold, define the unit vector $\hat z = z/\rho$.
In exact arithmetic, if the new row already lies in $\mathrm{span}(V)$ then $z=0$ and hence $\rho=0$,
so there is no meaningful new direction to append to the right singular subspace. In finite precision,
however, $\rho$ may be nonzero but extremely small due to roundoff; normalizing $z/\rho$ would then
amplify numerical noise and inject an essentially arbitrary direction into $V'$. We therefore treat
$\rho \le \mathrm{tol}$ as zero: the residual $z$ is dropped and no new direction is added, so Steps~2--3
below are applied without the last column of $K$ and without $\hat z^\top$, that is, with the
$(r+1)\times r$ core $\left[\begin{smallmatrix}\Sigma\\ p\end{smallmatrix}\right]$, its thin SVD, and
$V' = VV_K$. This keeps the update rank-preserving in this near-dependent case, keeps the columns of
$U'$ and $V'$ orthonormal, and improves numerical stability.

\paragraph{Step 2: small core SVD.}
Consider the block matrix
\[
K =
\begin{bmatrix}
\Sigma & 0 \\
p      & \rho
\end{bmatrix}
\in \mathbb{R}^{(r+1) \times (r+1)}.
\]
We compute a full SVD
$
K = U_K \Sigma' V_K^\top,
$
where $\Sigma' = \mathrm{diag}(\sigma'_1,\dots,\sigma'_{r+1})$.

\paragraph{Step 3: update singular vectors.}
The left singular vectors of the augmented matrix are obtained by ``lifting'' $U_K$ to the original dimension:
\[
U' =
\begin{bmatrix}
U & 0 \\
0 & 1
\end{bmatrix}
U_K \in \mathbb{R}^{(m+1) \times (r+1)},
\]
while the right singular vectors become
$
V' =
\begin{bmatrix}
V & \hat z^\top
\end{bmatrix}
V_K \in \mathbb{R}^{n \times (r+1)}.
$
The updated singular values are the diagonal entries of $\Sigma'$.

If the working rank $k$ is smaller than $r+1$, we truncate to the top $k$ singular values of $\Sigma'$ and the corresponding columns of $U'$ and $V'$. This yields a new truncated SVD of $A'$ at rank at most $k$.

The per-update cost is dominated by the SVD of the $(r+1) \times (r+1)$ core matrix $K$, which is $\mathcal{O}(r^3)$, plus the cost of the projections $x^\top V$ and the subsequent updates of $U$ and $V$, roughly $\mathcal{O}(nr^2 + mr^2)$. For small fixed $k$ this is essentially linear in $m+n$, in contrast to a full recomputation of the SVD, which costs $\mathcal{O}(mn \min\{m,n\})$.

\subsubsection{Column Append via Transpose}

Appending a new column $y \in \mathbb{R}^{m}$ is handled symmetrically by working on $A_t^\top$. 

A column append in $A_t$ corresponds to a row append in $A_t^\top$. In implementation, we conceptually form the SVD of the transpose, apply the row-append routine described above with the new row $y^\top$, and transpose back to obtain updated factors $(U',\Sigma',V')$ for the augmented matrix. Finally we truncate back to rank $k$ if needed.

This ``transpose trick'' reuses exactly the same algebra and implementation as the row case, and ensures that rows and columns are treated in a mathematically symmetric way.

\subsection{Incremental Rank-1 Entry Updates}
\label{subsec:rank1-entry-updates}

The third operation updates a single entry $(i,j)$ by a small increment $\delta$, as occurs when a user edits a rating or a single return is revised. Our algorithm follows the Brand-style update of~\cite{brand2006fast}. At the matrix level this is the rank-1 perturbation of the form \eqref{eq:rank1UpdateForm}.

Assume again that we have a rank-$r$ truncated SVD $A \approx U \Sigma V^\top$. Let
$
s_u^\top = U_{i,:} \in \mathbb{R}^{1 \times r}, %\qquad
s_v = V_{j,:}^\top \in \mathbb{R}^{r \times 1}
$
denote the $i$-th row of $U$ and the $j$-th row of $V$ (as a column vector) (i.e., the coordinates of $e_i$ and $e_j$ in the current latent bases), expressed in the current latent subspace. The update rule used in our implementation operates directly on the small $r \times r$ core. It forms the small core $K = \Sigma + \delta\, s_u s_v^\top$, computes its SVD $K = U_K \Sigma' V_K^\top$, and updates the factors via $U' = U U_K$, $V' = V V_K$, with new singular values $\Sigma'$.

This update is exact \emph{within} the current subspace in the following explicit sense. Let
$\widehat{A} := U\Sigma V^\top$ denote the current rank-$r$ approximation, and define the orthogonal
projectors onto the maintained left/right subspaces by $P_U := UU^\top$ and $P_V := VV^\top$
(since $U$ and $V$ have orthonormal columns). Writing $s_u = U^\top e_i$ and $s_v = V^\top e_j$,
the updated reconstruction produced by Steps 1--3 is
\[
U'\Sigma'(V')^\top
= U(U_K\Sigma'V_K^\top)V^\top
= UKV^\top
= U(\Sigma + \delta \, s_u s_v^\top) \, V^\top
= \widehat{A} + \delta\,(Us_u)(Vs_v)^\top.
\]
Using 
\[
Us_u = U(U^\top e_i) = (UU^\top) \, e_i = P_U e_i, \qquad Vs_v = (VV^\top) \, e_j = P_V e_j,
\]
we obtain
\begin{equation}
U'\Sigma'(V')^\top
= \widehat{A} + \delta\,(P_U e_i)(P_V e_j)^\top
= P_U(\widehat{A} + \delta \, e_ie_j^\top) \, P_V,
\label{eq:rank1-projection}
\end{equation}
which makes precise that we are applying the rank-1 perturbation and then projecting it back onto
$\mathrm{span}(U)\times\mathrm{span}(V)$. Any component of $\delta e_ie_j^\top$ outside the maintained
subspaces is intentionally discarded to keep rank fixed, matching the truncation viewpoint discussed
by Brand, who explicitly notes that rank-increasing components can be suppressed by ignoring the
part of an update that lies outside the current subspace~\cite{brand2006fast}. To this end, we proposed this rank suppression method.

The computational cost is dominated by the SVD of the $r\times r$ core $K$, i.e., $\mathcal{O}(r^3)$, plus the multiplications $UU_K$ and $VV_K$, which cost $\mathcal{O}(mr^2 + nr^2)$.

In practice we combine all three operations---row append, column append, and rank-1 perturbations---into a single incremental engine that supports heterogeneous update streams.

\subsection{Refreshes and Subspace Drift}
\label{subsec:refresh}

Because we truncate back to rank $k$ after each update, incremental operations inevitably introduce \emph{drift} between the maintained subspace and the ``optimal'' subspace one would obtain by recomputing a full SVD of the current matrix. Our framework makes this drift explicit and allows us to study different refresh strategies. This viewpoint aligns with classical subspace-tracking formulations, where drift is measured geometrically and stability over long streams is a central concern~\cite{moonen1992svdtrack}.

Alternative approaches handle drift differently: stochastic trackers such as GROUSE~\cite{balzano2010grouse} implicitly correct drift at every step via gradient projections but lack a mechanism to reset accumulated rounding errors, while Frequent Directions~\cite{liberty2013simple,ghashami2016frequent} maintains a provably near-optimal sketch but does not retain an explicit $U$ basis. Our refresh-based strategy occupies a middle ground, combining low per-update cost of deterministic Brand-style algebra with periodic batch recomputation to control accumulated error.

\subsubsection{Refresh Operation}

The \texttt{refresh} operation recomputes a truncated SVD of the current matrix:
$
A_t = U^{\text{full}} \, \Sigma^{\text{full}} \, (V^{\text{full}})^\top,
$
and then keeps only the top $k$ components by setting
$
U_t \gets U^{\text{full}}(:,1{:}k)
$,
$
\Sigma_t \gets \mathrm{diag}(\sigma^{\text{full}}_1,\dots,\sigma^{\text{full}}_k)
$, and
$
V_t \gets V^{\text{full}}(:,1{:}k).
$

We also maintain the scalar $N_t := \|A_t\|_F^2$ incrementally for use in the explained variance ratio (Section~\ref{subsec:metrics}), so that EVR can be evaluated without recomputing a full spectrum after each update.

In addition, the implementation can optionally store the current left singular vectors as a \emph{reference subspace},
$
U_{\text{ref}} \gets U_t,
$
which serves as the origin for measuring future subspace drift. When \texttt{refresh(store\_initial\_subspace = True)} is called, subsequent calls to \texttt{subspace\_principal\_angle} report the angle between $U_{\text{ref}}$ and the current $U_t$.

In the experiments we consider \emph{time-based refreshes} (e.g., every $R$ updates), the bound-triggered refresh of Section~\ref{subsec:bound-refresh}, and \emph{no-refresh} baselines, and we later discuss more adaptive strategies.

\subsubsection{Subspace Angles}

Given two matrices $Q_1, Q_2 \in \mathbb{R}^{m \times k}$ with orthonormal columns,
the principal angles between their column spaces are defined via the singular values
of $Q_1^\top Q_2$. By~\cite{bjorck1973numerical, golub1994canonical}, the singular values of $Q_1^\top Q_2$
are precisely the cosines of these principal angles, and since principal angles lie
in $[0, \tfrac{\pi}{2}]$ by definition, every singular value satisfies
$\sigma_\ell(Q_1^\top Q_2) \in [0, 1]$ for $\ell = 1, \dots, k$.
Let
\[
    Q_1^\top Q_2 = W \,\mathrm{diag}(\mu_1, \dots, \mu_k)\, Z^\top,
    \qquad 1 \ge \mu_1 \ge \dots \ge \mu_k \ge 0.
\]
Then the principal angles $\theta_\ell$ satisfy
$
    \cos\theta_\ell = \mu_\ell,\, \theta_\ell \in [0, \tfrac{\pi}{2}].
$

We focus on the largest principal angle, which corresponds to the smallest singular value:
\begin{equation}
\theta_{\max}(Q_1,Q_2) = \arccos \bigl(\min_\ell \mu_\ell \bigr)
= \arccos\bigl(\sigma_{\min}(Q_1^\top Q_2)\bigr).
\label{eq:theta-max}
\end{equation}
In our implementation we track two related quantities:
\begin{itemize}
    \item \textbf{Drift relative to a fixed reference subspace.}  
    With $Q_1 = U_{\text{ref}}$ and $Q_2 = U_t$, we define
    \begin{equation}
    \theta_{\text{ref}}(t)
    = \theta_{\max}(U_{\text{ref}}, U_t).
    \label{eq:angle-metrics-ref}
    \end{equation}
    This measures how far the current approximation has moved from the subspace established at the last refresh (or at initialization).

    \item \textbf{Angle to the optimal subspace.}  
    In the synthetic experiments we can also compute the full SVD of $A_t$ and obtain the optimal rank-$k$ subspace $U_{t}^{\text{opt}}$. We then define
    \begin{equation}
    \theta_{\text{opt}}(t)
    = \theta_{\max}(U_t^{\text{opt}}, U_t),
    \label{eq:angle-metrics-opt}
    \end{equation}
    which we record as \texttt{angle\_to\_opt}. This directly quantifies how close our incremental subspace is to the best possible one at time $t$.
\end{itemize}

These angles are central to our analysis of \emph{subspace stability} and serve as natural candidates for future \emph{adaptive refresh triggers} (e.g., refresh whenever $\theta_{\text{opt}}$ or $\theta_{\text{ref}}$ exceeds a threshold).

\subsection{Error Metrics and Computational Complexity}
\label{subsec:metrics}

For each state $A_t$ and its truncated SVD $U_t \Sigma_t V_t^\top$ we monitor several scalar metrics. All of them are implemented in the current codebase and reused across experiments.

\subsubsection{Frobenius Error and Optimal Baseline}

The basic approximation error is the Frobenius norm
$
E_{\text{inc}}(t)
= \| A_t - \widehat{A}_t \|_F,\quad
\widehat{A}_t := U_t \Sigma_t V_t^\top.
$

In synthetic settings we can also compute the \emph{optimal} rank-$k$ approximation using a full SVD of $A_t$. Let
$
A_t = U^{\text{full}} \, \Sigma^{\text{full}} \, (V^{\text{full}})^\top,
$
and define the best rank-$k$ approximation
$
A_t^{\text{opt}} = U^{\text{full}}(:,1{:}k)\,
\mathrm{diag}(\sigma^{\text{full}}_1,\dots,\sigma^{\text{full}}_k)\,
(V^{\text{full}}(:,1{:}k))^\top.
$
Having $\widehat{A}_t := U_t \Sigma_t V_t^\top$, we compute
\begin{equation}
E_{\text{inc}}(t)=\|A_t-\widehat{A}_t\|_F,
\qquad
E_{\text{opt}}(t)=\|A_t-A_t^{\text{opt}}\|_F.
\label{eq:frob-errors}
\end{equation}

Our metrics then include $\texttt{frob\_error}(t):=E_{\text{inc}}(t)$, $\texttt{frob\_opt}(t):=E_{\text{opt}}(t)$, $\texttt{frob\_gap}(t):=E_{\text{inc}}(t)-E_{\text{opt}}(t)$, and the error ratio
\begin{equation}
\texttt{frob\_ratio}(t)={E_{\text{inc}}(t)} \, / \, {E_{\text{opt}}(t)}.
\label{eq:frob-ratio}
\end{equation}

It is also worth noting that
$
E_{\text{inc}}(t) \ge E_{\text{opt}}(t).
$

\paragraph{Why $E_{\textnormal{inc}}(t)\ge E_{\textnormal{opt}}(t)$.}
By \eqref{eq:eckart-young-frob}, the best rank-$k$ approximation minimizes the Frobenius error over all rank-$k$ matrices. Since $\widehat{A}_t$ also has rank at most $k$, it follows that 
$
E_{\text{opt}}(t)=\|A_t-A_t^{\text{opt}}\|_F \le \|A_t-\widehat{A}_t\|_F = E_{\text{inc}}(t).
$

When $E_{\text{opt}}(t)$ is numerically zero we treat the ratio as undefined and record it as NaN. These quantities directly show how much accuracy we lose, if any, by updating incrementally instead of recomputing a full SVD at every step.

\subsubsection{Explained Variance Ratio}

We also report an explained variance ratio (EVR), defined as
\begin{equation}
\mathrm{EVR}(k;t)
=
\frac{\sum_{\ell=1}^k \sigma_{t,\ell}^2}
     {\sum_{\ell=1}^{r_t} \bigl(\sigma_{t,\ell}^{\text{full}}\bigr)^2}
=
\frac{\|\widehat{A}_t\|_F^2}{\|A_t\|_F^2},
\label{eq:evr}
\end{equation}
where $\widehat{A}_t = U_t \Sigma_t V_t^\top$ is the current maintained rank-$k$ approximation, $\sigma_{t,\ell}$ are its singular values, $\sigma_{t,\ell}^{\text{full}}$ are the singular values of the full matrix $A_t$, and $r_t$ is the numerical rank of $A_t$.

The singular-value form makes the interpretation explicit: the numerator sums the squared singular values retained in the current truncated model, while the denominator sums the squared singular values of the full matrix and therefore represents its total variance, equivalently its total Frobenius energy. Thus, EVR measures how much of the total matrix energy is retained by the rank-$k$ approximation currently maintained.

The equivalent Frobenius-norm form in \eqref{eq:evr} is especially convenient in practice. Since $U_t$ and $V_t$ have orthonormal columns, we have $\|\widehat{A}_t\|_F^2 = \|\Sigma_t\|_F^2 = \sum_{\ell=1}^k \sigma_{t,\ell}^2$, so the numerator is available directly from the maintained factors. For the denominator, we maintain the scalar $N_t := \|A_t\|_F^2$ incrementally. If a row $x^\top$ is appended, then $N_{t+1} = N_t + \|x\|^2$; if a column $y$ is appended, then $N_{t+1} = N_t + \|y\|^2$; and if a single entry is modified by $A_{t+1} = A_t + \delta e_i e_j^\top$, then only the $(i,j)$ entry changes, from $(A_t)_{ij}$ to $(A_t)_{ij} + \delta$, so
$N_{t+1} = N_t + ((A_t)_{ij}+\delta)^2 - (A_t)_{ij}^2 = N_t + 2\delta (A_t)_{ij} + \delta^2$.
Hence, $\|A_t\|_F^2$ can be tracked exactly without recomputing a full SVD.

In our synthetic experiments EVR mainly reflects how $k$ compares with the true rank, with or without refresh, rather than incremental drift: it stays above $0.99$ on the rank-1 and mixed streams, where $k$ is at least the true rank, and lies between $0.78$ and $0.83$ on the row and column streams of Section~\ref{subsec:structural}, where $k=4$ is below the true rank~6. In the finance application, EVR also provides a useful diagnostic of whether the chosen rank is sufficient to capture the dominant structure in the data.

\subsubsection{Runtime and Scaling}

Finally, the streaming simulator records the runtime spent per incremental update (\texttt{update\_time}) and, when applicable, the runtime spent computing the optimal baseline (\texttt{opt\_time}). We use these to study how latency scales with
\begin{itemize}
    \item the matrix size $(m_t, n_t)$, and
    \item the working rank $k$,
\end{itemize}
and to quantify the computational benefit of incremental updates versus full SVD recomputation.
All runtimes in this paper were measured on an Intel Core i7-13700H with one core in use (NumPy~2.4.4 with OpenBLAS~0.3.31 limited to one thread, process pinned to one performance core), as the median over five runs of the mean time per update; absolute timings are machine-dependent.
Full SVDs were timed in the same setup, as the mean time of a thin SVD (\texttt{numpy.linalg.svd} with \texttt{full\_matrices=False}) of a matrix of the same size, and are reported as multiples of the update time measured alongside.

The complexity of each incremental operation is dominated by SVDs of small $(k+1)\times(k+1)$ or $k\times k$ core matrices (for row/column and entry updates respectively), plus linear-in-size updates to the factors. In contrast, a full SVD of $A_t$ scales as $\mathcal{O}(m_t n_t \min\{m_t,n_t\})$, which rapidly becomes prohibitive in the streaming, large-scale regimes of interest.
Our implementation also stores the dense matrix $A_t$, which refreshes and the metrics need; each row or column append copies it, an $\mathcal{O}(m_tn_t)$ cost that is included in the measured update times (Section~\ref{subsec:finance-runtime}), whereas an entry update changes $A_t$ in place.

\subsection{A Posteriori Error Control Between Refreshes}
\label{subsec:error-control}

Between refreshes, evaluating $E_{\text{inc}}(t)$ requires forming $A_t - \widehat{A}_t$, and $E_{\text{opt}}(t)$ requires a full SVD.
We now show that the part of each update discarded by the engine has a norm that is available at negligible cost, and that these norms bound the incremental error and, on segments with a single update type, determine it.
Step $s$ of Algorithm~\ref{alg:incrementalsvd} maps $A_{s-1}$ to $A_s$; we write $E_s := A_s - \widehat{A}_s$, so that $E_{\text{inc}}(s) = \|E_s\|_F$, and let $t_0$ be the most recent refresh ($t_0 = 0$ for the initial truncated SVD).
Let $\widetilde{E}_{s-1}$ be $E_{s-1}$ padded with a zero row if step $s$ is a row append, padded with a zero column if it is a column append, and $\widetilde{E}_{s-1} = E_{s-1}$ for an entry update.
The part discarded at step $s$ and its size are
\[
D_s := E_s - \widetilde{E}_{s-1}, \qquad \eta_s := \|D_s\|_F .
\]

\begin{proposition}[A posteriori error control]
\label{prop:error-control}
Assume exact arithmetic. Let $U, \Sigma, V$ be the factors before step $s$, of rank $r \le k$, and $P_U = UU^\top$, $P_V = VV^\top$.
\begin{enumerate}
  \item[(a)] \emph{Discarded part.} For an entry update \eqref{eq:rank1UpdateForm},
  \[
  D_s = \delta\,\bigl(e_ie_j^\top - P_U e_ie_j^\top P_V\bigr), \qquad
  \eta_s = |\delta|\sqrt{1-ab}, \qquad a = \|U_{i,:}\|^2, \quad b = \|V_{j,:}\|^2 .
  \]
  For a row append with residual norm $\rho$ from \eqref{eq:row-residual}: if $\rho > \mathrm{tol}$, then $\eta_s = \sigma'_{k+1}$, the singular value of the core $K$ removed by the truncation ($\eta_s = 0$ if $r < k$); if $\rho \le \mathrm{tol}$, then $\eta_s = \rho$. In both cases $\eta_s \le \rho$. The same holds for a column append $y$, with $\rho = \|y - P_U y\|$ and $K$ the core of the transposed row append.
  \item[(b)] \emph{Cumulative bound.} For any sequence of updates and every $t \ge t_0$,
  \[
  E_{\textnormal{inc}}(t) \le \beta_t := E_{\textnormal{opt}}(t_0) + \sum_{s=t_0+1}^{t} \eta_s,
  \qquad
  E_{\textnormal{opt}}(t_0) = \bigl(N_{t_0} - \|\Sigma_{t_0}\|_F^2\bigr)^{1/2}.
  \]
  \item[(c)] \emph{Single-type segments.} If steps $t_0+1,\dots,t$ are all entry updates, or all row appends, or all column appends, with $\rho > \mathrm{tol}$ at every append, then
  \[
  E_{\textnormal{inc}}(t) = \bigl(N_t - \|\Sigma_t\|_F^2\bigr)^{1/2},
  \]
  and in the two append cases also $E_{\textnormal{inc}}(t)^2 = E_{\textnormal{opt}}(t_0)^2 + \sum_{s=t_0+1}^{t} \eta_s^2$.
\end{enumerate}
\end{proposition}

\begin{proof}
Write $\langle X, Y\rangle_F = \operatorname{tr}(X^\top Y)$, so that $\|X+Y\|_F^2 = \|X\|_F^2 + 2\langle X,Y\rangle_F + \|Y\|_F^2$.

(a) For an entry update, \eqref{eq:rank1-projection} gives $\widehat{A}_s = \widehat{A}_{s-1} + \delta\, P_U e_ie_j^\top P_V$, while $A_s = A_{s-1} + \delta\, e_ie_j^\top$; subtracting yields $D_s$.
Since $\|P_Ue_i\|^2 = e_i^\top P_U e_i = \|U^\top e_i\|^2 = a$, $\|P_Ve_j\|^2 = b$ and $\langle e_ie_j^\top, P_Ue_ie_j^\top P_V\rangle_F = (e_i^\top P_U e_i)(e_j^\top P_V e_j) = ab$, we get $\eta_s^2 = \delta^2(1 - 2ab + ab) = \delta^2(1-ab)$.
For a row append with $\rho > \mathrm{tol}$, let $\bar U = \left[\begin{smallmatrix} U & 0\\ 0 & 1\end{smallmatrix}\right]$ and $\bar V = [\,V \ \ \hat z^\top]$; both have orthonormal columns because $\hat z$ is a unit vector orthogonal to $\mathrm{span}(V)$.
Since $pV^\top + \rho\hat z = x^\top$, we have $\left[\begin{smallmatrix}\widehat{A}_{s-1}\\ x^\top\end{smallmatrix}\right] = \bar U K \bar V^\top$.
The engine keeps the leading $k$ singular triplets of $K$, i.e., $\widehat{A}_s = \bar U K_k \bar V^\top$ with $K_k$ the best rank-$k$ approximation of $K$, so
\[
D_s = \begin{bmatrix}A_{s-1}\\ x^\top\end{bmatrix} - \widehat{A}_s - \begin{bmatrix}E_{s-1}\\ 0\end{bmatrix}
= \begin{bmatrix}\widehat{A}_{s-1}\\ x^\top\end{bmatrix} - \widehat{A}_s
= \bar U (K - K_k) \bar V^\top
\]
and $\eta_s = \|K - K_k\|_F$. As $K$ has $r+1 \le k+1$ singular values, this is $\sigma'_{k+1}$ if $r = k$ and $0$ if $r < k$. Moreover $\sigma'_{r+1} = \sigma_{\min}(K) \le \|K e_{r+1}\| = \rho$, because the last column of $K$ is $\rho\, e_{r+1}$.
If $\rho \le \mathrm{tol}$, then $\bar U \left[\begin{smallmatrix}\Sigma\\ p\end{smallmatrix}\right] V^\top = \left[\begin{smallmatrix}\widehat{A}_{s-1}\\ x^\top - z\end{smallmatrix}\right]$ has only $r \le k$ singular triplets, so nothing is truncated, $\widehat{A}_s$ equals this matrix, and $D_s = \left[\begin{smallmatrix}0\\ z\end{smallmatrix}\right]$, i.e., $\eta_s = \rho$.
A column append is a row append for $A^\top$, with the roles of $U$ and $V$ exchanged.

(b) Since $E_s = \widetilde{E}_{s-1} + D_s$ and zero padding does not change the Frobenius norm, $E_{\text{inc}}(s) \le E_{\text{inc}}(s-1) + \eta_s$; summing over $s = t_0+1,\dots,t$ gives the bound.
At a refresh, $\widehat{A}_{t_0}$ is the rank-$k$ truncated SVD of $A_{t_0}$, so by \eqref{eq:eckart-young-frob} and the identities following it, $E_{\text{inc}}(t_0) = E_{\text{opt}}(t_0) = \bigl(\sum_{\ell>k}\sigma_\ell(A_{t_0})^2\bigr)^{1/2} = \bigl(N_{t_0} - \|\Sigma_{t_0}\|_F^2\bigr)^{1/2}$.

(c) Because $U_s$ and $V_s$ have orthonormal columns, $\|\widehat{A}_s\|_F = \|\Sigma_s\|_F$, so $N_s = \|\Sigma_s\|_F^2 + 2\langle E_s, \widehat{A}_s\rangle_F + \|E_s\|_F^2$; it suffices to show $\langle E_s, \widehat{A}_s\rangle_F = 0$ for $t_0 \le s \le t$.

\emph{Entry updates.} Here $U_s = U_{s-1}U_K$ and $V_s = V_{s-1}V_K$ with $U_K, V_K$ orthogonal, so $P_U$ and $P_V$ do not change on the segment. At the refresh, $U_{t_0}^\top A_{t_0} V_{t_0} = \Sigma_{t_0}$, i.e., $\widehat{A}_{t_0} = P_U A_{t_0} P_V$, and if $\widehat{A}_{s-1} = P_U A_{s-1} P_V$, then \eqref{eq:rank1-projection} gives $\widehat{A}_s = P_U (P_U A_{s-1} P_V + \delta e_ie_j^\top) P_V = P_U A_s P_V$. Since $X \mapsto P_U X P_V$ is an orthogonal projector for $\langle\cdot,\cdot\rangle_F$, $\langle E_s, \widehat{A}_s\rangle_F = \langle A_s - P_UA_sP_V,\, P_UA_sP_V\rangle_F = 0$.

\emph{Row appends.} We show by induction that $U_s^\top E_s = 0$, which gives $\langle E_s, \widehat{A}_s\rangle_F = \langle U_s^\top E_s, \Sigma_s V_s^\top\rangle_F = 0$. At the refresh, $U_{t_0}^\top A_{t_0} = \Sigma_{t_0} V_{t_0}^\top$, so $U_{t_0}^\top E_{t_0} = 0$. At step $s$, with the notation of (a), $U_s = \bar U\, U_K(:,1{:}k)$ (all columns of $U_K$ if $r < k$), and
\[
\bar U^\top \widetilde{E}_{s-1} = \begin{bmatrix} U_{s-1}^\top E_{s-1}\\ 0\end{bmatrix} = 0,
\qquad
U_K(:,1{:}k)^\top (K - K_k) = 0,
\]
the latter because the left singular vectors of $K$ are orthonormal. Hence $U_s^\top E_s = U_K(:,1{:}k)^\top \bigl(\bar U^\top\widetilde{E}_{s-1} + (K - K_k)\bar V^\top\bigr) = 0$. The first identity also gives $\langle \widetilde{E}_{s-1}, D_s\rangle_F = \langle \bar U^\top \widetilde{E}_{s-1} \bar V, K - K_k\rangle_F = 0$, so $\|E_s\|_F^2 = \|E_{s-1}\|_F^2 + \eta_s^2$; summing from $s = t_0+1$ to $t$ and using $\|E_{t_0}\|_F = E_{\text{opt}}(t_0)$ gives the second identity.

\emph{Column appends.} Apply the row case to $A^\top$; the invariant becomes $E_s V_s = 0$.
\end{proof}

\paragraph{Sharpness.}
The bound in (b) is attained: for $A_0 = \mathrm{diag}(1, \tfrac12)$, $k = 1$ and repeated updates of entry $(2,2)$ by $\delta > 0$, every step has $a = b = 0$ and $E_{\text{inc}}(t) = \beta_t$, whereas $\bigl(E_{\text{opt}}(0)^2 + \sum_{s \le t} \eta_s^2\bigr)^{1/2} < \beta_t$; hence the root-sum-of-squares identity in (c) does not extend to entry updates.
Nor does (c) extend to mixed segments: updating entry $(1,2)$ of the same $A_0$ by $\delta = 1$ and then appending the row $(1, 1)$ gives $E_{\text{inc}}(2) \approx 0.859$, while $\bigl(N_2 - \|\Sigma_2\|_F^2\bigr)^{1/2} = \bigl(E_{\text{opt}}(0)^2 + \eta_1^2 + \eta_2^2\bigr)^{1/2} \approx 1.278$.

\paragraph{Use as a refresh trigger.}
Every quantity in Proposition~\ref{prop:error-control} is available without a full SVD: $\eta_s$ costs $\mathcal{O}(k)$ for an entry update and is a by-product of the update for an append ($\sigma'_{k+1}$ or $\rho$), $N_t$ is updated in $\mathcal{O}(1)$ as described after \eqref{eq:evr}, $\|\Sigma_t\|_F^2$ costs $\mathcal{O}(k)$, and $E_{\text{opt}}(t_0)$ is the norm of the singular values discarded at the last refresh.
The oracle ratio in the error-based rule \eqref{eq:error-refresh-rule} can therefore be replaced by the computable ratio $\beta_t / E_{\text{opt}}(t_0)$, or by $E_{\text{inc}}(t) / E_{\text{opt}}(t_0)$ on single-type segments.
For append-only streams this replacement is safe: $A_{t_0}$ is a submatrix of $A_t$, so $E_{\text{opt}}(t) \ge E_{\text{opt}}(t_0)$ (restricting a rank-$k$ approximation of $A_t$ to that submatrix gives one of $A_{t_0}$), and refreshing as soon as the computable ratio exceeds $\gamma$ keeps $\texttt{frob\_ratio}(t) \le \gamma$ at every step, provided $E_{\text{opt}}(t_0) > 0$.

\medskip

Taken together, this framework provides a clean separation between (i) the \emph{algebraic core} of the incremental SVD updates (Sections~\ref{subsec:row-col-updates}--\ref{subsec:rank1-entry-updates}), (ii) the \emph{refresh} mechanism and subspace drift diagnostics (Section~\ref{subsec:refresh}), (iii) the \emph{metrics and complexity} used to evaluate accuracy and stability (Section~\ref{subsec:metrics}), and (iv) the computable a posteriori bound on the error between refreshes (Section~\ref{subsec:error-control}). The next section turns these ingredients into refresh and rank-selection policies; Section~\ref{sec:synthetic} instantiates the framework on controlled synthetic streams, and Section~\ref{sec:finance} applies the same engine to financial data.

\section{Refresh and Rank-Selection Policies}
\label{sec:policies}
The row/column append updates in Section~\ref{subsec:row-col-updates} are exact up to rank-$k$
truncation at each step, while the rank-1 entry update in Section~\ref{subsec:rank1-entry-updates}
performs an exact \emph{in-subspace} correction at fixed rank. Over long streams two effects emerge:
(i) approximation error can accumulate, and (ii) the ``intrinsic'' rank and dominant subspace of the underlying matrix may change.
In practice one therefore needs policies that decide
\emph{when} to recompute a fresh SVD (``refresh'') and
\emph{how large} the maintained rank $k$ should be.

Related streaming approaches that update low-rank structure via alternative factorizations likewise emphasize the need for occasional re-alignment to maintain accuracy and numerical stability over long streams~\cite{brust2025svdtype}.

This section formalizes a design space of such policies, expressed
directly in terms of the metrics defined in Section~\ref{subsec:metrics}:
reconstruction error $\texttt{frob\_error}$, optimal error
$\texttt{frob\_opt}$, error ratio $\texttt{frob\_ratio}$,
subspace angles $\texttt{angle}$ and $\texttt{angle\_to\_opt}$,
explained-variance ratio $\mathrm{EVR}(k)$, and runtime measures
$\texttt{update\_time}$ and $\texttt{opt\_time}$.
In the synthetic and application experiments (Sections~\ref{sec:synthetic}--\ref{sec:finance}) we instantiate several of these policies and
compare their behavior.

\subsection{Periodic Refresh}
\label{sec:periodic-refresh}

The simplest policy is to refresh on a fixed schedule.
Let $t$ denote the number of updates processed so far (rank-1, row, or column
updates), and fix an integer refresh period $T_{\mathrm{refresh}} \ge 1$.
A \emph{periodic refresh} policy is defined by
\begin{equation}
\text{refresh at time } t
\quad \Longleftrightarrow \quad
t \equiv 0 \pmod{T_{\mathrm{refresh}}}.
\label{eq:periodic-refresh-rule}
\end{equation}
Operationally, a refresh recomputes a truncated SVD of the current matrix $A_t$ and replaces $(U_t,\Sigma_t,V_t)$ by the top-$k$ singular triplets returned by the full SVD routine. In practice, we call the same \texttt{refresh()} routine that is used for
sanity checks in our implementation.

This policy exposes an intuitive trade-off:
\begin{itemize}
  \item Between refreshes the algorithm uses only cheap incremental updates
  with per-step cost $\texttt{update\_time}$.
  \item At refresh steps it pays a one-off cost $\texttt{opt\_time}$ to
  recompute the SVD.
\end{itemize}
On long streams, small $T_{\mathrm{refresh}}$ (frequent refreshes) improves
numerical fidelity at the expense of more full SVD calls,
whereas large $T_{\mathrm{refresh}}$ reduces the number of full SVDs but
allows more drift to accumulate.
Section~\ref{sec:synthetic} evaluates this trade-off in a controlled
rank-1 streaming setting.

\subsection{Error-Based Refresh}
\label{sec:error-refresh}

While periodic refresh uses time as the only signal, our framework also
provides direct error measurements via the baseline metrics.

Whenever a full SVD of $A_t$ is available, we can evaluate the optimal rank-$k$ baseline error $E_{\text{opt}}(t)$ from \eqref{eq:frob-errors}; by \eqref{eq:eckart-young-frob}, this is the minimum Frobenius error over all rank-$k$ approximations.

We then monitor the relative degradation through the ratio $\texttt{frob\_ratio}(t)$ defined in \eqref{eq:frob-ratio}, where $\widehat{A}_t$ denotes the current incremental approximation and $A_t^{\text{opt}}$ the best rank-$k$ truncation from the full SVD.

An \emph{error-based} refresh policy triggers whenever the incremental error
drifts too far from the optimal baseline:
\begin{equation}
\text{refresh at time } t
\quad \Longleftrightarrow \quad
\texttt{frob\_ratio}(t) > \gamma,
\label{eq:error-refresh-rule}
\end{equation}
for some tolerance $\gamma > 1$ (e.g.\ $\gamma \in [1.1,1.2]$).
To avoid excessive chattering, it is natural to combine this with a minimum
spacing constraint $T_{\min}$:
we only allow refreshes for $t$ satisfying
$t - t_{\text{last}} \ge T_{\min}$, where $t_{\text{last}}$
is the time of the most recent refresh.

In synthetic experiments we can compute $\texttt{frob\_opt}(t)$ exactly
at every step and therefore treat error-based refresh as an ``oracle''
policy.
In real systems, recomputing a full SVD at each step defeats the purpose of
being incremental.
There, error-based strategies can still be approximated by:
(i) computing $\texttt{frob\_opt}$ only on a coarse schedule
(e.g.\ once per hour), or
(ii) using cheaper surrogates for the optimal error, such as tracking the
residual norm of new rows or monitoring changes in the spectrum of
$\Sigma_t$ over time.
Section~\ref{subsec:error-control} gives a computable replacement for the oracle ratio: $\beta_t/E_{\text{opt}}(t_0)$, where $t_0$ is the last refresh, needs no full SVD between refreshes and is an upper bound on $\texttt{frob\_ratio}(t)$ on append-only streams.
On such streams, refreshing as soon as $\beta_t/E_{\text{opt}}(t_0) > \gamma$ keeps $\texttt{frob\_ratio}(t) \le \gamma$ at every step, provided $E_{\text{opt}}(t_0) > 0$; with entry updates $E_{\text{opt}}(t)$ can fall below $E_{\text{opt}}(t_0)$, and the guarantee is lost.
The rule we test in Sections~\ref{subsec:bound-refresh} and~\ref{sec:finance}, the bound-triggered refresh \eqref{eq:bound-refresh-rule}, instead compares $\beta_t$ with $\|A_t\|_F$: it refreshes when $\beta_t/\|A_t\|_F > \tau$, which keeps $E_{\text{inc}}(t) \le \tau\|A_t\|_F$ after every step as long as $\tau$ exceeds $E_{\text{opt}}(t)/\|A_t\|_F$, but does not bound $\texttt{frob\_ratio}(t)$.

\subsection{Angle-Based Refresh}
\label{sec:angle-refresh}

Reconstruction error focuses on \emph{values} in the matrix.
For applications such as latent-factor modeling, the geometry of the
subspace spanned by the top singular vectors is equally important.
Our metrics therefore include two principal-angle based measures:
\begin{itemize}
  \item $\texttt{angle}(t)$:
  the principal angle between the current subspace
  $\mathcal{U}_t = \operatorname{span}(U_t)$
  and the \emph{initial} subspace $\mathcal{U}_0$.
  This captures cumulative drift since initialization.
  \item $\texttt{angle\_to\_opt}(t)$:
  the principal angle between $\mathcal{U}_t$ and the
  optimal rank-$k$ subspace
  $\mathcal{U}_t^{\text{opt}} = \operatorname{span}(U_t^{\text{opt}})$
  obtained from the full SVD of $A_t$.
  This measures how far the incremental subspace is from the best
  possible one at time $t$.
\end{itemize}

Our angular diagnostics are based on the largest principal angle in \eqref{eq:theta-max}; in particular, we monitor the quantities in \eqref{eq:angle-metrics-ref}--\eqref{eq:angle-metrics-opt}.
An \emph{angle-based} refresh policy declares the current subspace outdated
once the misalignment exceeds a user-specified threshold:
\[
  \text{refresh at time } t
  \quad \Longleftrightarrow \quad
  \texttt{angle\_to\_opt}(t) > \theta_{\max},
\]
for some $\theta_{\max} \in (0, \tfrac{\pi}{2}]$.
Compared to pure error-based triggers, angle-based policies are more directly
tied to subspace tracking and are less sensitive to global rescaling of
the matrix.

As for error-based refresh, computing $\texttt{angle\_to\_opt}(t)$ requires
access to the optimal subspace and is therefore mainly feasible in synthetic
or diagnostic settings.
However, the same idea can be approximated in practice without a full SVD
by using local angular signals.
For example:
\begin{itemize}
  \item \emph{New-row angle:}
  when a new row $x_t^\top$ arrives, we can decompose it as shown in \eqref{eq:row-residual}, then the
  angle between $x_t$ and its projection $p_t$,
  or equivalently the ratio $\|z_t\| \, / \, \|x_t\|$, serves as a measure of how
  ``novel'' this observation is relative to the current subspace.
  Large values indicate concept drift.
  
  \item \emph{Subspace-change per update:}
  the angle between successive subspaces
  $\operatorname{span}(U_{t-1})$ and $\operatorname{span}(U_t)$ is cheap
  to compute and can signal instability if it becomes consistently large.
\end{itemize}
In Section~\ref{sec:discussion} we interpret such angular signals as a form
of concept-drift detector for dynamic matrices.

\subsection{Adaptive Rank Selection}
\label{sec:adaptive-rank}

So far we have treated the target rank $k$ as fixed. In many applications, however, the intrinsic dimensionality of the data changes over time. Our framework therefore supports adaptive rank selection through the explained-variance ratio in \eqref{eq:evr}, the residual quantities from \eqref{eq:row-residual}, and principal-angle diagnostics derived from \eqref{eq:theta-max}.

\paragraph{EVR-based rule.}
Given singular values $\sigma_1(t) \ge \dots \ge \sigma_r(t)$ of $A_t$,
the explained-variance ratio defined in \eqref{eq:evr} gives
a natural rule: choose the smallest $k$ such that
$\mathrm{EVR}(k; t) \ge \tau_{\mathrm{EVR}}$ for some target level
(e.g.\ $\tau_{\mathrm{EVR}} = 0.9$ or $0.95$).
In practice we enforce bounds $k_{\min} \le k_t \le k_{\max}$ and update
$k_t$ only when $\mathrm{EVR}(k_t; t)$ drifts significantly away from the
target to avoid oscillations:
\begin{equation}
k_t
=
\min\left\{
k \in [k_{\min},k_{\max}]
:
\mathrm{EVR}(k;t) \ge \tau_{\mathrm{EVR}}
\right\}.
\label{eq:evr-rank-rule}
\end{equation}

\paragraph{Novelty- and residual-based rules.}
EVR captures global retained energy, but in streaming settings it is also useful to react to \emph{local novelty} in newly arriving data. For a newly appended row $x_t^\top$, the decomposition in \eqref{eq:row-residual} provides a natural signal: the residual component $z_t$ measures the part of the new observation that lies outside the current right-singular subspace. Consequently, the relative residual
$\|z_t\| \, / \, \|x_t\|$
quantifies how poorly the current rank-$k$ model explains the incoming row.

A simple adaptive rule is therefore to increase the maintained rank when this novelty signal exceeds a prescribed threshold. Concretely, one can use
\begin{equation}
\|z_t\| > \nu \,\|x_t\|
\quad \Longrightarrow \quad
k_{t+1} \gets \min(k_t+1,\,k_{\max}),
\label{eq:novelty-rank-rule}
\end{equation}
where $\nu \in (0,1)$ controls the sensitivity of the trigger and $k_{\max}$ is a user-specified upper bound on the admissible rank.

This rule complements the EVR-based criterion: while EVR indicates whether the maintained approximation retains enough total matrix energy overall, the residual test reacts to directions that are new relative to the current subspace, even when the global energy ratio remains high. In this sense, EVR captures global adequacy, whereas the residual-based rule captures local subspace novelty.

\paragraph{Interaction with refresh.}
Adaptive rank selection and refresh policies are complementary.
A refresh step may:
(i) recompute a rank-$k_t$ SVD with the current adaptive rank, or
(ii) temporarily compute a higher-rank SVD to reassess $\mathrm{EVR}$
and update $k_t$.

\subsection{Summary of Policy Design Space}
\label{sec:policy-summary}

Table~\ref{tab:policy-summary} summarizes the main policy families discussed
above.
All policies are expressed in terms of the metrics produced by our
implementation (Section~\ref{subsec:metrics}), and can be plugged into the same
incremental SVD engine without modifying its core updates.
The experimental sections compare two of these policies, periodic and
bound-triggered refresh, in terms of accuracy, subspace drift, and runtime.

\begin{table}[t]
  \centering
  \caption{Overview of refresh and rank-selection policies.
  Each policy is defined in terms of the metrics from
  Section~\ref{subsec:metrics} or the bound of Section~\ref{subsec:error-control}}
  \label{tab:policy-summary}
  \vspace{0.3em}
  \begin{tabularx}{\linewidth}{l X X}
    \toprule
    Policy type & Signal(s) used & Pros / cons \\
    \midrule
    Periodic refresh
      & Time index $t$; period $T_{\mathrm{refresh}}$
      & Simple, predictable cost; ignores actual error or drift. \\
    Error-based refresh
      & \texttt{frob\_ratio}($t$), optional $T_{\min}$
      & Direct control of approximation error; requires (approximate) baseline. \\
    Bound-triggered refresh
      & $\beta_t/\|A_t\|_F$ (Proposition~\ref{prop:error-control}); tolerance $\tau$
      & $\mathcal{O}(1)$ test, no baseline; certifies $E_{\text{inc}}(t)\le\tau\|A_t\|_F$ if $\tau$ exceeds the floor $E_{\text{opt}}(t)/\|A_t\|_F$, but not \texttt{frob\_ratio}. \\
    Angle-based refresh
      & \texttt{angle}($t$), \texttt{angle\_to\_opt}($t$), new-row angles
      & Targets subspace drift; aligns with factor interpretations; requires angular estimates. \\
    Adaptive rank
      & $\mathrm{EVR}(k;t)$, residual norms, novelty angles
      & Adjusts model complexity to data; introduces extra hyperparameters ($\tau_{\mathrm{EVR}}, k_{\min}, k_{\max}$). \\
    \bottomrule
  \end{tabularx}
\end{table}

\section{Synthetic Experiments}
\label{sec:synthetic}

We now evaluate the incremental SVD framework on controlled synthetic
streams generated from known low-rank models. The goals are to measure:
(i) accuracy of the incremental approximation relative to the optimal
rank-$k$ baseline; (ii) subspace drift and its control via refresh
policies; (iii) sensitivity to the truncation rank $k$; (iv) behavior
under structural growth; and (v) runtime scaling.

All experiments use the ground-truth generators in
\texttt{experiments/streams.py}, the metrics in \texttt{isvd.metrics}, and the
incremental engine \texttt{isvd.IncrementalSVD} of the released Python package
\texttt{isvd} (see Code availability). At each logged step we compute a full SVD of
the current dense matrix, yielding the optimal rank-$k$ error
(\texttt{frob\_opt}) and optimal subspace for comparison.

\subsection{Rank-1 Streaming: Tracking and Drift}
\label{subsec:rank1}

We begin with a $50\times40$ matrix of true rank~5 and perform $10{,}000$
random rank-1 updates of size $\delta\sim\mathcal{N}(0,0.05^2)$.
We compare no refresh with the periodic refresh policy \eqref{eq:periodic-refresh-rule} using $T_{\mathrm{refresh}}=1000$.

\begin{figure}[!htbp]
  \centering
  \begin{subfigure}[t]{\panelfigwidth}
    \centering
    \includegraphics[width=\linewidth]{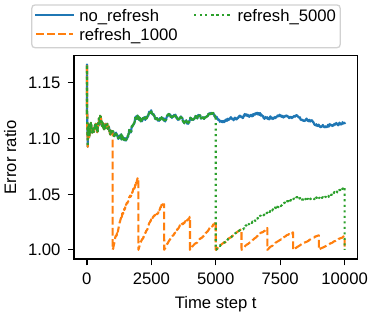}
    \caption{Error ratio $\texttt{frob\_ratio}(t)$ from \eqref{eq:frob-ratio}. Refresh keeps the ratio near~1}
    \label{fig:rank1_error_ratio}
  \end{subfigure}
  \hspace{\hspacegapsize\textwidth}
  \begin{subfigure}[t]{\panelfigwidth}
    \centering
    \includegraphics[width=\linewidth]{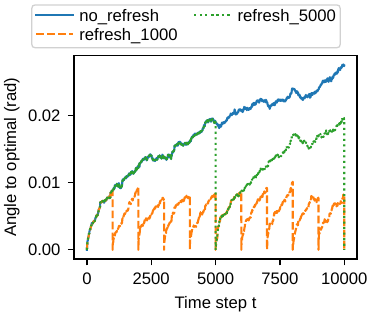}
    \caption{Angle to the optimal subspace, $\texttt{angle\_to\_opt}(t)$ from \eqref{eq:angle-metrics-opt}. Refresh repeatedly realigns the basis}
    \label{fig:rank1_angle_opt}
  \end{subfigure}
  \caption{Rank-1 streaming: drift diagnostics under different refresh policies: (a) relative Frobenius error ratio; (b) angle to the optimal rank-$5$ subspace. Solid: no refresh; dashed: refresh every 1000 steps; dotted: refresh every 5000 steps}
  \label{fig:rank1_diagnostics}
\end{figure}

Figure~\ref{fig:rank1_error_ratio} shows the error ratio $\texttt{frob\_ratio}(t)$ from \eqref{eq:frob-ratio}.
Without refresh, the incremental approximation drifts to roughly $1.12$,
while periodic refresh keeps the ratio extremely close to~1.
Figure~\ref{fig:rank1_angle_opt} shows the angle to the optimal rank-$5$ subspace, namely $\texttt{angle\_to\_opt}(t)$ from \eqref{eq:angle-metrics-opt}: between refreshes the angle grows slowly, then resets
to nearly zero, forming the characteristic ``sawtooth'' pattern. This confirms that incremental updates accumulate moderate drift, but
periodic refresh tightly controls the alignment to the true subspace.

\subsection{Effect of Truncation Rank}
\label{subsec:k-effect}

\begin{figure}[!htbp]
  \centering
  \includegraphics[width=\singlefigwidth]{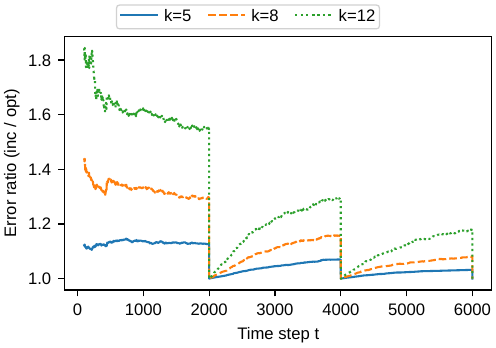}
  \caption{Effect of $k$ on error ratio, refresh every 2000 steps ($k=5$ solid, $k=8$ dashed, $k=12$ dotted)}
  \label{fig:k_error_ratio}
\end{figure}

We next study how the truncation rank $k$ interacts with incremental drift.
We run a rank-1 stream of the same kind as in Section~\ref{subsec:rank1} ($50\times40$, true rank~5), with 6000 updates and a refresh every 2000 steps, for
$k \in \{5,8,12\}$ and plot the error ratio $\texttt{frob\_ratio}(t)$ from \eqref{eq:frob-ratio} for $t\ge100$ in Figure~\ref{fig:k_error_ratio}.
The vertical resets at $t=2000$ and $4000$ are the refreshes.

At first glance the figure looks counter-intuitive: the curves for larger
$k$ sit \emph{above} the curves for smaller $k$ (mean ratio 1.06, 1.16 and
1.30 for $k=5$, $8$ and $12$, over $t\ge100$), so the error ratio increases
with $k$. This does not mean that higher rank makes the
approximation worse. In fact, the absolute Frobenius errors in \eqref{eq:frob-errors} decrease with $k$, as expected. To interpret the ratio, let $E_{\text{inc}}^{(k)}(t)$ and $E_{\text{opt}}^{(k)}(t)$ denote the quantities in \eqref{eq:frob-errors} for working rank $k$, and write $E_{\text{inc}}^{(k)}(t)=E_{\text{opt}}^{(k)}(t)+\delta_k(t)$ with $\delta_k(t)\ge 0$. Then $\texttt{frob\_ratio}(t,k)=1+\delta_k(t) \, / \, E_{\text{opt}}^{(k)}(t)$.

As $k$ increases from the true latent rank (here $5$), the optimal rank-$k$ error $E_{\text{opt}}^{(k)}(t)$ falls (mean 2.32, 1.99 and 1.63 for $k=5$, $8$ and $12$), and so does the incremental error $E_{\text{inc}}^{(k)}(t)$, but more slowly (mean 2.44, 2.24 and 2.02). The incremental gap $\delta_k(t)$ therefore grows with $k$ (mean 0.12, 0.24 and 0.39), and the fraction $\delta_k(t) \, / \, E_{\text{opt}}^{(k)}(t)$ grows with it, because its numerator grows while its denominator shrinks.

There is also a modeling aspect.  For $k$ larger than the true rank, the
additional singular directions mainly capture high-frequency noise.  The
full SVD re-estimates these noise directions from scratch at each step and
therefore achieves the best possible rank-$k$ fit.  The incremental method,
by contrast, tries to track these rapidly changing, low-energy directions
using small Brand-style updates; this is inherently harder and contributes
to the growing gap $\delta_k(t)$.  The angle to the optimal rank-$k$ subspace,
$\texttt{angle\_to\_opt}$ from \eqref{eq:angle-metrics-opt}, shows the same:
its mean is 0.011 radians for $k=5$, but 1.02 and 1.24 radians for $k=8$
and $12$, where it reaches $\tfrac{\pi}{2}$.  From a practical point of
view, this suggests that increasing $k$ well beyond the intrinsic rank
offers little benefit: it reduces the optimal error, but mainly on noisy
components that are difficult for any incremental scheme to track, and it
raises both computational cost and the relative error ratio.

\begin{figure}[!htbp]
  \centering
  \begin{subfigure}[t]{\panelfigwidth}
    \centering
    \includegraphics[width=\linewidth]{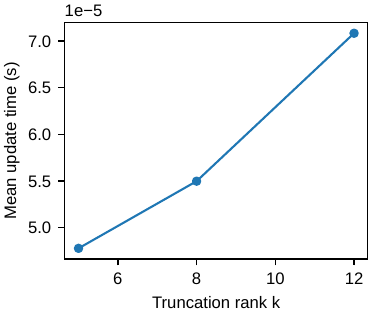}
    \caption{Mean incremental update time vs.\ $k$, rank-1 stream on a $50\times40$ matrix. Runtime grows by a factor of 1.5 from $k=5$ to $k=12$}
    \label{fig:k_time}
  \end{subfigure}
  \hspace{\hspacegapsize\textwidth}
  \begin{subfigure}[t]{\panelfigwidth}
    \centering
    \includegraphics[width=\linewidth]{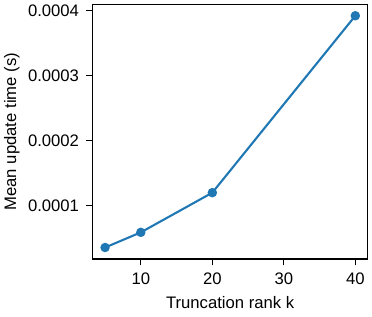}
    \caption{Scaling with truncation rank $k$, rank-1 updates on a $500\times400$ matrix. Runtime grows faster than linearly in $k$}
    \label{fig:scaling_rank}
  \end{subfigure}
  \caption{Runtime scaling with truncation rank $k$ in two synthetic settings (median over five runs, one core). Runtime increases with $k$ in both; the growth is clearly faster than linear only on the larger matrix (b)}
  \label{fig:rank_runtime_comparison}
\end{figure}

The incremental update time, shown in
Figure~\ref{fig:k_time}, increases with~$k$: the mean is $4.8$, $5.5$ and
$7.1\cdot 10^{-5}$\,s per update for $k=5$, $8$ and $12$, a factor of 1.5
from $k=5$ to $k=12$.  On this small $50\times 40$ matrix the cost therefore
grows much more slowly than the $\mathcal{O}(k^3)$ term of the core SVD alone
would suggest (a factor of about 14).
This highlights the practical speed--accuracy trade-off.

\subsection{Structural Growth: Rows and Columns}
\label{subsec:structural}

We next consider expanding matrices.
New rows or columns are sampled from a rank-6 model with additive noise.
We track a rank-4 approximation and refresh every 50 additions.

\begin{figure}[!htbp]
  \centering
  
  \begin{subfigure}[t]{\panelfigwidth}
    \centering
    \includegraphics[width=\linewidth]{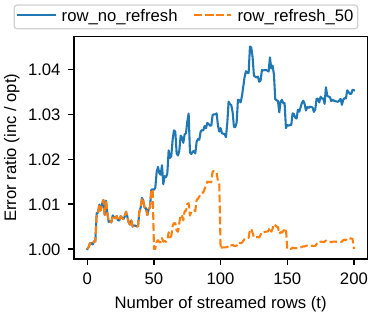}
    \caption{Row streaming}
    \label{fig:row_error_ratio}
  \end{subfigure}
  \hspace{\hspacegapsize\textwidth}
  \begin{subfigure}[t]{\panelfigwidth}
    \centering
    \includegraphics[width=\linewidth]{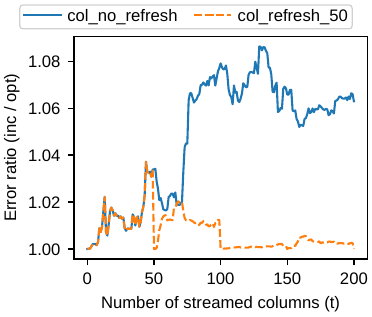}
    \caption{Column streaming}
    \label{fig:col_error_ratio}
  \end{subfigure}

  \caption{Structural growth: error ratio under no refresh (solid) vs.\ refresh every 50 steps (dashed) for (a) rows and (b) columns.
  Refresh maintains near-optimal error}
  \label{fig:rowcol_error_ratio}
\end{figure}

Figure~\ref{fig:rowcol_error_ratio} shows the error ratio for both row and
column streaming.
In both cases, no-refresh gradually drifts away from the optimal error,
while periodic refresh maintains $\texttt{frob\_ratio}(t)$ from \eqref{eq:frob-ratio} close to 1.0--1.02.
This demonstrates that the incremental engine remains robust under
changing matrix dimensions, and that refresh policies remain effective
in the structurally growing setting.

\subsection{Mixed Streaming}
\label{subsec:mixed}

\begin{figure}[!htbp]
  \centering
  \includegraphics[width=\singlefigwidth]{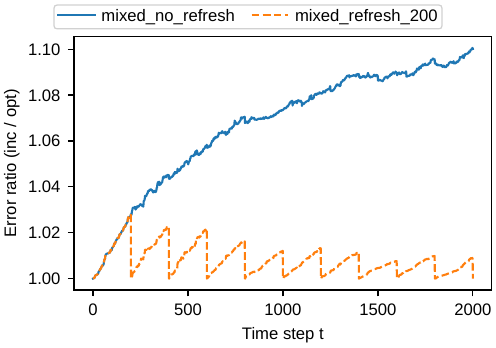}
  \caption{Mixed streaming: error ratio without refresh (solid) and with refresh every 200 steps (dashed).
  Without refresh, drift compounds across update types;
  refresh stabilizes the approximation}
  \label{fig:mixed_error_ratio}
\end{figure}

We combine rank-1 perturbations with row additions (one new row every 50 steps).
Figure~\ref{fig:mixed_error_ratio} shows the error ratio $\texttt{frob\_ratio}(t)$ from \eqref{eq:frob-ratio}.
We see that the unrefreshed method
accumulates more than 10\% error relative to the optimal baseline,
while refresh every 200 steps keeps the ratio very close to~1.
This is consistent with the more complex drift patterns induced by the
heterogeneous stream, and highlights the importance of periodic
realignment.

\subsection{A Posteriori Bound and Bound-Triggered Refresh}
\label{subsec:bound-refresh}

We return to the rank-1 stream of Section~\ref{subsec:rank1} ($50\times40$, true rank~5, $k=5$, $10{,}000$ updates, $\delta\sim\mathcal{N}(0,0.05^2)$) and record the bound $\beta_t$ of Proposition~\ref{prop:error-control}(b), which the engine updates in $\mathcal{O}(1)$ per step.
Each policy runs on five seeds; seed~0 is the stream of Section~\ref{subsec:rank1}.

\begin{figure}[!htbp]
  \centering
  \includegraphics[width=\singlefigwidth]{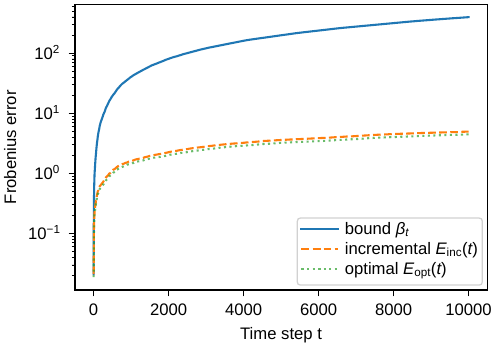}
  \caption{Rank-1 streaming without refresh (seed~0): true error $E_{\text{inc}}(t)$ (dashed), a posteriori bound $\beta_t$ of Proposition~\ref{prop:error-control}(b) (solid), and optimal error $E_{\text{opt}}(t)$ (dotted), log scale}
  \label{fig:rank1_error_bound}
\end{figure}

\paragraph{Tightness.}
Figure~\ref{fig:rank1_error_bound} shows $E_{\text{inc}}(t)$, $\beta_t$ and $E_{\text{opt}}(t)$ without refresh.
The bound holds at every step of all 30 runs reported here, up to roundoff: the largest value of $E_{\text{inc}}(t)-\beta_t$ is $7.7\times10^{-14}$.
It is also loose.
Across the seeds, $\beta_t/E_{\text{inc}}(t)$ is 25.0--25.5 at $t=1000$ and 77.9--81.4 at $t=10{,}000$.
On seed~0, $\beta_t$ grows linearly (40.4 at $t=1000$, 396.5 at $t=10{,}000$) and $E_{\text{inc}}$ roughly like $\sqrt{t}$ (1.60 and 4.96).
The triangle inequality behind $\beta_t$ adds up the $\eta_s$ and ignores any cancellation between the discarded parts $D_s$ of different steps.
A refresh resets $\beta_t$ to $E_{\text{opt}}(t_0)$, and with $T_{\mathrm{refresh}}=1000$ the median of $\beta_t/E_{\text{inc}}(t)$ drops to 7.6 (Table~\ref{tab:bound-refresh}).
On this stream, which has only entry updates, Proposition~\ref{prop:error-control}(c) gives $E_{\text{inc}}(t)$ exactly, so $\beta_t$ matters mostly for mixed streams.

\paragraph{Bound-triggered refresh.}
We refresh after step $t$ when the relative bound exceeds a tolerance $\tau$,
\begin{equation}
\beta_t \,/\, \|A_t\|_F > \tau,
\qquad \|A_t\|_F = N_t^{1/2}.
\label{eq:bound-refresh-rule}
\end{equation}
The test costs $\mathcal{O}(1)$ per step.
It guarantees $E_{\text{inc}}(t)\le\tau\|A_t\|_F$ after every step as long as $\tau$ exceeds the floor $E_{\text{opt}}(t)/\|A_t\|_F$, because a refresh brings $\beta_t$ down to $E_{\text{opt}}(t)$ and no lower.
Table~\ref{tab:bound-refresh} compares the rule with periodic refresh.
Each extra refresh buys accuracy: the mean error ratio falls from 1.118 without refresh to 1.004 with $\tau=0.05$, and the mean angle falls with it.
At an equal number of full SVDs, $\tau=0.4$ and $T_{\mathrm{refresh}}=1000$ are practically tied (mean ratio 1.025 and 1.026).
This is expected here, since $\beta_t$ grows at a nearly constant rate on this stream and the rule therefore refreshes at nearly even intervals.
Smaller $\tau$ pushes the refreshes towards the end of the stream, because the floor rises from $0$ to 0.042--0.052 at $t=10{,}000$.
With $\tau=0.05$, on average 375 of the 448 refreshes fall in the second half.
On one seed the floor ends above $\tau$ (0.0517), and the rule then refreshes at every step, 1322 times in total.
No policy does much about the worst step.
Every run reaches its maximum ratio before its first refresh, in the first segment, which starts from $E_{\text{opt}}(0)=0$, so the mean of the maxima is 1.176--1.180 for every policy.
After the first refresh the maximum is 1.048--1.079 in every run.
The rule bounds the error relative to $\|A_t\|_F$, not relative to $E_{\text{opt}}(t)$.
For $\tau\ge0.1$ the guaranteed level is also far above the actual error: in every run, even without refresh, $E_{\text{inc}}(t)/\|A_t\|_F$ stays below 0.058.
In practice $\tau$ sets the refresh budget rather than the accuracy.
On append-only streams the trigger $\beta_t/E_{\text{opt}}(t_0)$ of Section~\ref{subsec:error-control} bounds \texttt{frob\_ratio} itself; on this stream of entry updates it does not, because $E_{\text{opt}}(t)$ can fall below $E_{\text{opt}}(t_0)$ (Section~\ref{sec:error-refresh}).

\begin{table}[!htbp]
  \centering
  \caption{Rank-1 streaming (Section~\ref{subsec:rank1} setup): periodic refresh vs.\ the bound-triggered rule (refresh when $\beta_t/\|A_t\|_F>\tau$), mean over seeds 0--4.
  SVDs: number of refreshes (excluding the initial SVD), with the range over seeds; $t>T/2$: refreshes in the second half.
  Ratio and angle statistics are over $t=1,\dots,10{,}000$, angles in $10^{-3}$\,rad.
  Last column: median of $\beta_t/E_{\text{inc}}(t)$}
  \label{tab:bound-refresh}
  \vspace{0.3em}
  \footnotesize
  \setlength{\tabcolsep}{2.5pt}
  \begin{tabular}{l r r c c c c c c}
    \toprule
    & \multicolumn{2}{c}{SVDs} & \multicolumn{3}{c}{\texttt{frob\_ratio}} & \multicolumn{2}{c}{\texttt{angle\_to\_opt}} & \\
    \cmidrule(lr){2-3}\cmidrule(lr){4-6}\cmidrule(lr){7-8}
    Policy & total & $t>T/2$ & mean & max & max, $t>T/2$ & mean & max & $\beta_t/E_{\text{inc}}$ \\
    \midrule
    No refresh & 0 & 0 & 1.118 & 1.180 & 1.122 & 19.6 & 29.4 & 56.2 \\
    Periodic, $T_{\mathrm{refresh}}=1000$ & 10 & 5 & 1.026 & 1.180 & 1.020 & 6.2 & 10.9 & 7.6 \\
    Bound, $\tau=0.4$ & 10.4 (10--12) & 5.6 & 1.025 & 1.180 & 1.021 & 5.9 & 10.1 & 7.1 \\
    Bound, $\tau=0.2$ & 22.8 (21--27) & 12.4 & 1.014 & 1.177 & 1.010 & 4.1 & 7.5 & 3.7 \\
    Bound, $\tau=0.1$ & 57.4 (52--70) & 32.6 & 1.008 & 1.177 & 1.005 & 2.7 & 5.6 & 2.1 \\
    Bound, $\tau=0.05$ & 447.8 (214--1322) & 374.8 & 1.004 & 1.176 & 1.001 & 1.4 & 4.1 & 1.2 \\
    \bottomrule
  \end{tabular}
\end{table}

\subsection{Efficiency and Scaling}
\label{subsec:scaling}

For fixed $k$, the rank-1 update time hardly depends on the matrix size at these dimensions: for $k=5$ it is $4.8\cdot10^{-5}$\,s on the $50\times40$ matrix of Section~\ref{subsec:k-effect} and $3.6\cdot10^{-5}$\,s on the $500\times400$ matrix of Figure~\ref{fig:scaling_rank}, which has 100 times as many entries, so at these sizes the $\mathcal{O}((m+n)k^2)$ updates of the factors are a small part of the cost. Figure~\ref{fig:scaling_rank} ($500\times400$, 300 rank-1 updates per $k$) shows faster-than-linear growth in $k$: the mean update time is $3.6$, $5.9$, $12.0$ and $39.2\cdot10^{-5}$\,s for $k=5$, $10$, $20$ and $40$, so doubling $k$ from 20 to 40 multiplies it by 3.3, less than the factor 8 of a pure $\mathcal{O}(k^3)$ cost. This is consistent with the complexity discussion in Section~\ref{subsec:metrics} and with the small-core update structure in Sections~\ref{subsec:row-col-updates}--\ref{subsec:rank1-entry-updates}, where both the core SVD and the updates of the factors grow with $k$. The absolute runtimes remain below a millisecond for the problem sizes considered. A full SVD costs far more: it takes about 6 times as long as one rank-1 update with $k=5$ on the $50\times40$ matrix, and about 1400 times as long on the $500\times400$ matrix.

\medskip

Across all scenarios, incremental SVD updates remain accurate,
stable, and computationally efficient.
Periodic and bound-triggered refresh both bound drift relative to the
optimal subspace; the bound-triggered rule also certifies
$E_{\text{inc}}(t)\le\tau\|A_t\|_F$ after every step as long as $\tau$
exceeds the floor $E_{\text{opt}}(t)/\|A_t\|_F$, and the engine adapts robustly to structural growth
in matrix dimensions.

\section{Financial Factor and Covariance Modeling}
\label{sec:finance}

In our application we treat incremental SVD as an online equity
\emph{factor model}: each new trading day adds a row of returns, and we
update a low-rank approximation of the return matrix in real time.
From a quant perspective this corresponds to maintaining PCA-style
systematic risk factors and a low-rank covariance matrix used for
portfolio construction and risk management.

The standard approaches for online covariance estimation in finance include rolling-window PCA, dynamic conditional correlation models (DCC-GARCH)~\cite{engle2002dcc}, shrinkage estimators~\cite{ledoit2004wellconditioned}, and thresholded factor models such as POET~\cite{fan2013large}. DCC-GARCH directly models time-varying correlations but scales quadratically in the number of assets; recent work unifies the DCC and shrinkage lines into a single DCC-NL framework augmented with intraday OHLC information, scaling to universes of $N \geq 1000$ assets for Markowitz portfolio selection~\cite{denard2022largedynamic}. Shrinkage and POET produce well-conditioned covariance matrices but are typically applied in batch on a rolling window. A concurrent streaming alternative is to combine a factor model with a knowledge-based sketch trained on historical data, yielding direct covariance recovery from sketched streaming returns and demonstrating competitive portfolio-allocation performance against POET~\cite{tan2024keef}; this is closely related to our approach but relies on offline-trained sketches rather than incremental SVD updates. Incremental SVD offers a complementary approach: it updates the factor structure itself one observation at a time, yielding an explicit low-rank covariance at every step without re-estimating on the full window.

\subsection{Data and Preprocessing}
\label{subsec:finance-data}

We construct a multi-asset universe of $N=60$ liquid US-listed ETFs,
covering broad equity indices, sector and style funds, credit and rates
ETFs, and a few commodity proxies.\footnote{Data are downloaded from
Yahoo Finance using a small Python script; the exact ticker list is
available with the code.}
For each asset we download daily adjusted close prices and align them on
the intersection of trading days.  After removing days with any missing
prices we obtain a panel of about $T \approx 1{,}640$ trading days,
i.e.\ roughly six to seven years of history.

From prices $P_{t,i}$ we build log-returns $r_{t,i}=\log(P_{t,i} \, / \, P_{t-1,i})$, stacked into a returns matrix $R\in\mathbb{R}^{T\times N}$ (time $\times$ assets).

Before running SVD we perform \emph{causal} de-meaning (Option~A) to avoid look-ahead bias in the online setting.  For each asset $i$ we maintain an expanding mean
\[
  \textstyle \mu_{t,i} := \frac{1}{t}\sum_{s=1}^{t} r_{s,i},
\qquad
  R(t,i) := r_{t,i} - \mu_{t,i},
\]
so that each return is centered using only information available up to
(and including) time $t$.  In a streaming implementation this is updated
online via the standard recurrence
\[
  \mu_t = \mu_{t-1} + \tfrac{1}{t}\bigl(r_t - \mu_{t-1}\bigr),
\qquad
  R_t = r_t - \mu_t,
\]
where $\mu_t, r_t \in \mathbb{R}^N$ denote the cross-sectional mean and
return vectors at day $t$, and $R_t$ is the centered row appended to the
matrix.  For readability we keep the notation $R$ for this causally
centered returns matrix throughout the remainder of Section~\ref{sec:finance}.

This is standard in factor modeling: the sample covariance of centered
returns is then
\begin{equation}
\Sigma_{\mathrm{full}}(t)=\tfrac{1}{t-1}  R_{1:t}^\top \, R_{1:t},
\label{eq:sample-covariance}
\end{equation}
and the
right singular vectors of $R_{1:t}$ coincide with the eigenvectors of
$\Sigma_{\mathrm{full}}(t)$.  We deliberately do \emph{not} standardize by individual
volatilities; the low-rank model works directly with covariance in
return units, which is what a risk engine needs.

A common alternative in live risk systems is to center returns using an
\emph{exponentially weighted} mean (Option~B), which discounts older data
and adapts more quickly under regime changes:
\[
  \mu^{\text{EW}}_t = (1-\alpha) \, \mu^{\text{EW}}_{t-1} + \alpha r_t,
\qquad
  R_t = r_t - \mu^{\text{EW}}_t,
\]
with a forgetting factor $\alpha \in (0,1)$ (often specified via a
half-life).  Like Option~A, this update is fully causal and can be applied
online with $\mathcal{O}(N)$ state.  In this paper we use the expanding-mean
centering (Option~A) for the main results because it is parameter-free;
Option~B is a drop-in replacement when faster adaptation is desired.

Note that using a full-sample mean (centering by $T^{-1} \sum_{s=1}^T r_{s,i}$) would introduce look-ahead information; the causal centering above avoids this while preserving the online nature of the experiment.

Although we use Yahoo data for convenience, the pipeline only requires a
clean time$\times$asset returns matrix.  The same incremental SVD code
would plug directly into institutional data sources such as WRDS/CRSP by
replacing the data loading step.

\subsection{Online Factor Model via Incremental SVD}
\label{subsec:finance-factor}

We interpret the columns of $V_t \in \mathbb{R}^{N \times k}$ in the SVD
$R_{1:t} \approx U_t \Sigma_t V_t^\top$ as $k$ orthogonalized
\emph{factor portfolios} (PCA factors), and the rows of
$U_t \Sigma_t$ as time series of factor returns.
To mimic a live deployment we split the history into
\begin{itemize}
  \item an \emph{initial window} of $T_0 = 250$ trading days
        (roughly one year) used to compute a dense SVD and initialize
        the incremental model, and
  \item a streaming phase where we append one new row of returns per day
        and update the truncated SVD.
\end{itemize}
We focus first on a fixed truncation rank $k=5$, which is a very common
dimension for PCA-style equity factor models (one ``market mode'' plus
a handful of sector/curve/style directions).  For the streaming phase, we
compare three refresh policies, instances of \eqref{eq:periodic-refresh-rule}:
\begin{itemize}
  \item \textbf{no\_refresh}: pure incremental SVD, never recomputing a
        full SVD after $T_0$;
  \item \textbf{refresh\_100}: periodic refresh every 100 trading days
        (roughly quarterly);
  \item \textbf{refresh\_20}: aggressive refresh every 20 trading days
        (approximately monthly).
\end{itemize}
Every few days we also compute the \emph{optimal} rank-$k$ SVD of the
\emph{same causally centered} matrix $R_{1:t}$ and use it as an oracle
baseline.  This is not something a production system would do, but it allows us to quantify how well the incremental SVD tracks the best possible factor subspace.

\begin{figure}[!htbp]
  \centering
  \includegraphics[width=\singlefigwidth]{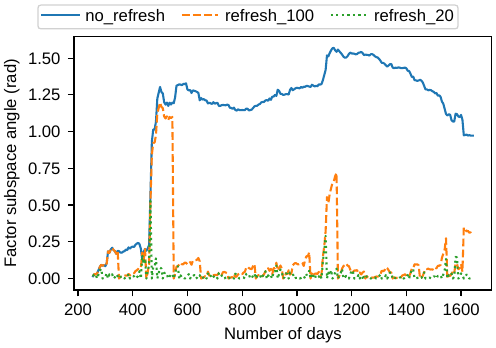}
  \caption{Finance application, $k=5$:
  angle between incremental and oracle factor subspaces over time
  without refresh (solid) and with refresh every 100 days (dashed) and every 20 days (dotted)}
  \label{fig:finance-angle-factor}
\end{figure}

Figure~\ref{fig:finance-angle-factor} shows the principal angle between
the incremental factor subspace and the oracle subspace introduced in \eqref{eq:angle-metrics-opt} as a function of
time.  With no refresh, the angle drifts steadily, approaching values
near $\tfrac{\pi}{2}$ radians: the incremental factors eventually become almost
orthogonal to the current optimal PCA factors.  Periodic refresh every
100 days keeps the angle bounded and produces visible ``sawtooth''
patterns: drift accumulates between refreshes and is reset when we
recompute the SVD.  Refreshing every 20 days keeps the factor subspace
close to the oracle most of the time (median angle 0.013\,rad, under one
degree), with two short spikes after abrupt market moves: 0.53\,rad on
day 465 (late April 2020, after the COVID crash) and 0.29\,rad on day
1105 (November 2022), each removed by the next refresh.  From a quant perspective this means
that the incremental model continues to capture the dominant systematic
risk directions (market, sectors, rates, etc.) even as the return
structure evolves.

\subsection{Covariance and Portfolio Risk Approximation}
\label{subsec:finance-risk}

Given the causally centered returns matrix $R_{1:t}$ we define the sample covariance as in \eqref{eq:sample-covariance}.
Let $R_{1:t} = U_t \Sigma_t V_t^\top$ be the full SVD. Then
\begin{equation}
\Sigma_{\mathrm{full}}(t)
= \tfrac{1}{t-1}(U_t \Sigma_t V_t^\top)^\top (U_t \Sigma_t V_t^\top)
= \tfrac{1}{t-1}V_t \, \Sigma_t^2 \, V_t^\top.
\label{eq:cov-from-svd}
\end{equation}
Thus, $R_{1:t}^\top \, R_{1:t}$ has eigendecomposition $V_t \Sigma_t^2 V_t^\top$:
its eigenvectors are the right singular vectors of $R_{1:t}$ and its eigenvalues are $\sigma_i^2(t)$.
Therefore, 
$\Sigma_{\mathrm{full}}(t)$
has eigenvalues $\sigma_i^2(t) \, / \, (t-1)$.

Let $V_{k,\mathrm{full}}$ be the first $k$ right singular vectors with singular
values $\sigma_{1:k}(t)$.  The oracle low-rank covariance is
\begin{align}
\widehat{\Sigma}_{\mathrm{full},k}(t)
&=
\tfrac{1}{t-1} \,
V_{k,\mathrm{full}} \cdot
{\mathrm{diag}(\sigma_1^2(t),\dots,\sigma_k^2(t))} \cdot
V_{k,\mathrm{full}}^\top,
\nonumber\\
\widehat{\Sigma}_{\mathrm{inc},k}(t)
&=
\tfrac{1}{t-1} \, V_{k,\mathrm{inc}} \cdot
{\mathrm{diag}(\widetilde{\sigma}_1^2(t),\dots,\widetilde{\sigma}_k^2 (t))} \cdot
V_{k,\mathrm{inc}}^\top,
\label{eq:low-rank-covariances}
\end{align}
where $V_{k,\mathrm{inc}}$ and $\widetilde{\sigma}_i(t)$ come from the
incremental SVD.  This is the object that would feed into a Markowitz
or risk-parity engine.

\begin{figure}[!htbp]
  \centering
  \begin{subfigure}[t]{\panelfigwidth}
    \centering
    \includegraphics[width=\linewidth]{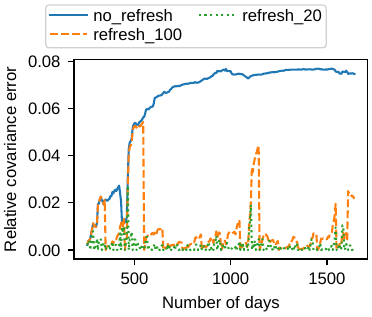}
    \caption{Covariance relative error vs.\ time}
    \label{fig:finance_cov_error}
  \end{subfigure}
  \hspace{\hspacegapsize\textwidth}
  \begin{subfigure}[t]{\panelfigwidth}
    \centering
    \includegraphics[width=\linewidth]{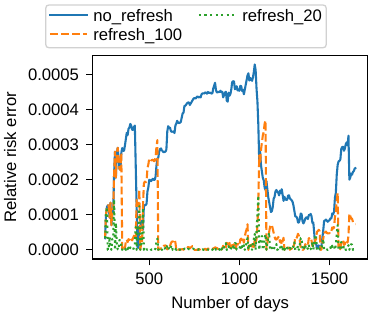}
    
    \caption{Equal-weight portfolio risk error vs.\ time}
    \label{fig:finance_risk_eq}
  \end{subfigure}

  \caption{Finance application, $k=5$: covariance and risk approximation
  without refresh (solid) and with refresh every 100 days (dashed) and every 20 days (dotted). Errors are measured against the
  oracle rank-$k$ model obtained from a full SVD of $R_{1:t}$}
  \label{fig:finance-cov-risk}
\end{figure}

\paragraph{Covariance error.}
We measure covariance quality via a relative Frobenius error
\begin{equation}
\mathrm{err}_{\Sigma}(t)
=
\frac{\|\widehat{\Sigma}_{\mathrm{inc},k}(t)-\widehat{\Sigma}_{\mathrm{full},k}(t)\|_F}
     {\|\widehat{\Sigma}_{\mathrm{full},k}(t)\|_F}.
\label{eq:cov-error}
\end{equation}
Figure~\ref{fig:finance-cov-risk}(a) shows
$\mathrm{err}_\Sigma(t)$ for $k=5$.  Without refresh the relative error
gradually climbs to around $7$–$8\%$ by the end of the sample.  In other
words, the low-rank covariance implied by the incremental factors
eventually deviates meaningfully from the best rank-5 covariance.  A
quarterly refresh (every 100 days) caps the error at a few percent and
regularly brings it back down after each recomputation.  Monthly
refreshes keep the error at the sub-percent level (median $0.12\%$) except
in the same two short episodes as the angle, with peaks of $2.7\%$ on day
465 and $1.9\%$ on day 1105, and the final relative error is about
$6\cdot 10^{-4}$ (also for $T_{\mathrm{refresh}}=30$ in Table~\ref{tab:finance-bound-refresh}).  This
matches the subspace picture: tighter alignment of factor loadings
translates directly into a more accurate covariance surface.

A notable feature in Figure~\ref{fig:finance-cov-risk}(a) is the sharp
rise in covariance error between days 460 and 500 (late April to
mid-June 2020): without refresh it grows from $1.4\%$ to $5.4\%$.
It follows the 2020 COVID market crash (days 420--442 of our sample,
20 February to 23 March 2020), in which cross-asset correlations and
volatilities jumped abruptly over a few trading days.
This kind of regime break is effectively a \emph{structural shock} to
the covariance matrix, and neither an incremental model nor an oracle
rank-$k$ approximation can ``predict'' it ex ante.
What we see instead is that all low-rank models incur a transient bump
in error when the new regime arrives.
The key difference is in recovery: with periodic refresh (especially
\texttt{refresh\_20}) the incremental SVD quickly re-aligns to the new
eigenstructure, bringing covariance errors back to the sub-percent
range, whereas \texttt{no\_refresh} continues to drift as the market
environment evolves after the shock.

\paragraph{Portfolio risk.}
For trading and risk management, what matters is the impact on
\emph{portfolio} risk, not the Frobenius norm of the covariance matrix.
Given a portfolio weight vector $w \in \mathbb{R}^N$ we define
the oracle and incremental volatilities
\begin{equation}
\sigma_{\mathrm{full}}(t)=\sqrt{\smash[b]{w^\top \widehat{\Sigma}_{\mathrm{full},k}(t) w}},
\qquad
\sigma_{\mathrm{inc}}(t)=\sqrt{\smash[b]{w^\top \widehat{\Sigma}_{\mathrm{inc},k}(t) w}},
\label{eq:portfolio-vol}
\end{equation}
and track the relative error
\begin{equation}
\mathrm{err}_{\mathrm{risk}}(t)
=
\frac{|\sigma_{\mathrm{inc}}(t)-\sigma_{\mathrm{full}}(t)|}
     {\sigma_{\mathrm{full}}(t)}.
\label{eq:risk-error}
\end{equation}
Figure~\ref{fig:finance-cov-risk}(b) reports this quantity for an
equal-weight long-only portfolio $w_{\text{eq}} = \mathbf{1}/N$.  The
picture is striking: even when the covariance Frobenius error for
\texttt{no\_refresh} reaches several percent, the corresponding relative
error in portfolio volatility stays at or below $5.3\cdot 10^{-4}$ (roughly
$0.05\%$).  Periodic refreshes reduce this by another order of magnitude
or more; for example, with $k=5$ and refresh every 30 or 20 trading
days, final relative risk errors are in the $10^{-6}$--$10^{-5}$ range
for the equal-weight portfolio.  At every evaluation we also draw 10
random long--short portfolios, weight vectors with independent standard
normal entries scaled to unit norm, and record the mean of
$\mathrm{err}_{\mathrm{risk}}$ over them.  These errors are about two
orders of magnitude larger than for $w_{\text{eq}}$ but follow the same
ordering: on the final day the mean is $3.8\cdot 10^{-2}$ without
refresh, $1.8\cdot 10^{-2}$ with refresh every 100 days and
$3.9\cdot 10^{-4}$ with refresh every 20 days.

The same episode appears in Figure~\ref{fig:finance-cov-risk}(b) as
a rise in risk error after day 460 for the two refreshed runs.
Realized volatilities and correlations moved so quickly that factors
calibrated on the pre-crash regime misestimate risk for a while.
Even in this stress window (days 450--600) the largest relative error for
the equal-weight portfolio is $1.2$--$3.4\cdot 10^{-4}$ across the three
policies.  How long it lasts depends on the next refresh:
\texttt{refresh\_20} removes its peak ($1.2\cdot 10^{-4}$ on day 465) five
days later, whereas \texttt{refresh\_100} keeps errors of
$1.2$--$3.1\cdot 10^{-4}$ until its refresh on day 550 (late August 2020).
From a trading and risk perspective this is exactly what one would
expect: extreme tail events are not forecast by a Gaussian factor model,
but an online low-rank engine can adapt quickly once the new regime is
observed.

\begin{figure}[!htbp]
  \centering
  \begin{subfigure}[t]{\panelfigwidth}
    \centering
    \includegraphics[width=\linewidth]{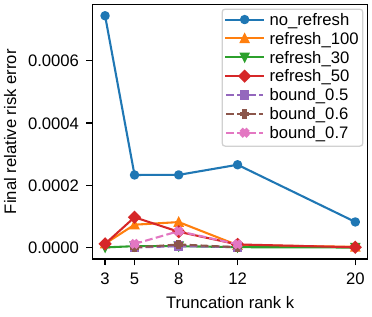}
    \caption{Final relative risk error of the equal-weight portfolio
    as a function of truncation rank $k$ and refresh policy.  Periodic
    refresh plus moderate $k$ yields very small errors.  Solid: no refresh
    (circles) and refresh every 100 (upward triangles), 30 (downward
    triangles) and 50 days (diamonds). Dashed: the
    bound-triggered rule \eqref{eq:bound-refresh-rule},
    $\tau=0.5$ (squares), $0.6$ (plus signs) and $0.7$ (crosses), run for $k\in\{5,8,12\}$}
    \label{fig:finance-rank-risk}
  \end{subfigure}
  \hspace{\hspacegapsize\textwidth}
  \begin{subfigure}[t]{\panelfigwidth}
    \centering
    \includegraphics[width=\linewidth]{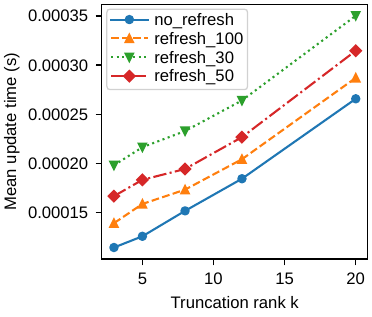}
    \caption{Mean incremental update time vs.\ rank $k$ without refresh
    (solid, circles) and with refresh every 100 (dashed, upward triangles),
    30 (dotted, downward triangles) and 50 days (dash-dotted, diamonds). Runtime increases approximately linearly with
    $k$, while more frequent refresh shifts the curves upward by a
    modest amount}
    \label{fig:finance-runtime}
  \end{subfigure}
  \caption{Finance application, rank--refresh trade-off across truncation ranks and refresh cadences: (a) accuracy of the equal-weight risk estimate; (b) per-update runtime cost}
  \label{fig:finance-rank-risk-runtime}
\end{figure}

\paragraph{Rank--refresh trade-off.}
To explore the interaction between truncation rank and refresh cadence
we run a small grid over $k \in \{3,5,8,12,20\}$ and
$\texttt{refresh\_every} \in \{\texttt{None}, 30, 50, 100\}$.
For each configuration we record the final-day angle between factor
subspaces, covariance error, portfolio risk error, and mean update time.
The equal-weight risk errors are summarized in
Figure~\ref{fig:finance-rank-risk}.  A few patterns emerge:

\begin{itemize}
  \item With \texttt{no\_refresh}, increasing $k$ from 5 to 20 lowers the
        final covariance error from $7.5\cdot 10^{-2}$ to $7.7\cdot 10^{-3}$
        and the equal-weight risk error from $2.3\cdot 10^{-4}$ to
        $8.2\cdot 10^{-5}$, the latter not monotonically ($2.7\cdot 10^{-4}$
        at $k=12$), while the factor drift stays large (final angles of
        $0.52$--$0.97$ radians for $k=5$--$20$).  At $k=3$ the final angle is
        only $0.084$ radians, but the risk error is the largest in the grid
        ($7.4\cdot 10^{-4}$).
  \item For periodic refreshes, moderate ranks $k=5$ or $k=8$ combined
        with a monthly or 30-day refresh already deliver extremely small
        risk errors (down to $\mathcal{O}(10^{-6})$) and covariance
        errors well below $0.1\%$.
  \item Very frequent refreshes yield diminishing accuracy returns while
        increasing runtime; in this ETF universe the sweet spot is
        around $k \in [5, 8]$ with refresh every 20--50 days.
\end{itemize}

From a quant point of view, this means that a compact low-rank factor
model with a handful of factors and occasional re-baselining is
sufficient to keep risk estimates for diversified portfolios extremely
accurate.

\paragraph{Bound-triggered refresh.}
We also run the rule \eqref{eq:bound-refresh-rule} of Section~\ref{subsec:bound-refresh} for $k\in\{5,8,12\}$ with $\tau\in\{0.5,0.6,0.7\}$ (Table~\ref{tab:finance-bound-refresh}, dashed curves in Figure~\ref{fig:finance-rank-risk}).
The return matrix is far from low rank, so the floor $E_{\text{opt}}(t)/\|A_t\|_F$ is high: its maximum over the stream is 0.451, 0.362 and 0.289 for $k=5,8,12$.
All three tolerances lie above it, but the guarantee $E_{\text{inc}}(t)\le\tau\|A_t\|_F$ is then weak, and $\tau$ acts as a refresh budget.
The margin between $\tau$ and the floor grows with $k$, so at fixed $\tau$ the rule refreshes less often for larger $k$: 196, 90 and 52 full SVDs at $\tau=0.5$ and 63, 41 and 28 at $\tau=0.7$, against 13--46 for the periodic policies.
The final-day errors follow the number of days since the last full SVD more closely than the policy.
For $k=5$, $\tau=0.5$ and $\tau=0.6$ both refresh two days before the end and give the smallest errors in the table (angle $1.1\cdot10^{-3}$\,rad, $\mathrm{err}_\Sigma=8.6\cdot10^{-5}$, equal-weight risk error $3.5\cdot10^{-7}$), at the cost of 196 and 96 full SVDs.
For $k=8$, $\tau=0.5$ and $T_{\mathrm{refresh}}=30$ last refresh on the same day and give identical final values.
At similar SVD counts the two rules are close when their last refreshes are about equally recent ($k=12$: $\tau=0.7$ with 28 SVDs and $T_{\mathrm{refresh}}=50$ with 27, final angles 0.056 and 0.057\,rad).
At $k=8$, $\tau=0.7$ with 41 SVDs ends worse than $T_{\mathrm{refresh}}=30$ with 46 (0.061 vs.\ 0.015\,rad), because its last refresh is 44 days old rather than 13.
A single final-day snapshot therefore does not separate the two rules on this stream.
What the bound-triggered rule adds is the certificate $E_{\text{inc}}(t)\le\tau\|A_t\|_F$ after every step, not a smaller final error.

\begin{table}[!htbp]
  \centering
  \caption{Finance application, final day ($t=1{,}643$, causal de-meaning): periodic refresh vs.\ the bound-triggered rule \eqref{eq:bound-refresh-rule}.
  SVDs: number of refreshes in the 1{,}393 streaming days (excluding the initial SVD at $T_0=250$); days: days since the last full SVD.
  Angle: largest principal angle \eqref{eq:angle-metrics-opt} in $10^{-3}$\,rad; $\mathrm{err}_\Sigma$ \eqref{eq:cov-error} in $10^{-4}$; $\mathrm{err}_{\mathrm{risk}}$ \eqref{eq:risk-error} for $w_{\text{eq}}$ in $10^{-6}$, and its mean over 10 random unit-norm Gaussian weight vectors in $10^{-4}$}
  \label{tab:finance-bound-refresh}
  \vspace{0.3em}
  \small
  \setlength{\tabcolsep}{4pt}
  \begin{tabular}{r l r r r r r r}
    \toprule
    & & & & & & \multicolumn{2}{c}{$\mathrm{err}_{\mathrm{risk}}$} \\
    \cmidrule(lr){7-8}
    $k$ & Policy & SVDs & days & angle & $\mathrm{err}_\Sigma$ & $w_{\text{eq}}$ & random \\
    \midrule
    5 & No refresh & 0 & 1393 & 973 & 746 & 233 & 379 \\
     & Periodic, $T_{\mathrm{refresh}}=30$ & 46 & 13 & 7.66 & 6.05 & 4.65 & 3.93 \\
     & Periodic, $T_{\mathrm{refresh}}=50$ & 27 & 43 & 397 & 284 & 97.4 & 238 \\
     & Periodic, $T_{\mathrm{refresh}}=100$ & 13 & 93 & 304 & 221 & 73.5 & 181 \\
     & Bound, $\tau=0.5$ & 196 & 2 & 1.10 & 0.859 & 0.349 & 0.697 \\
     & Bound, $\tau=0.6$ & 96 & 2 & 1.10 & 0.859 & 0.349 & 0.697 \\
     & Bound, $\tau=0.7$ & 63 & 23 & 16.1 & 12.5 & 11.9 & 8.99 \\
    \midrule
    8 & No refresh & 0 & 1393 & 766 & 246 & 233 & 58.5 \\
     & Periodic, $T_{\mathrm{refresh}}=30$ & 46 & 13 & 15.0 & 6.28 & 4.71 & 2.38 \\
     & Periodic, $T_{\mathrm{refresh}}=50$ & 27 & 43 & 60.3 & 23.4 & 50.9 & 10.7 \\
     & Periodic, $T_{\mathrm{refresh}}=100$ & 13 & 93 & 110 & 43.1 & 82.0 & 20.9 \\
     & Bound, $\tau=0.5$ & 90 & 13 & 15.0 & 6.28 & 4.71 & 2.38 \\
     & Bound, $\tau=0.6$ & 57 & 20 & 20.0 & 7.97 & 10.1 & 3.07 \\
     & Bound, $\tau=0.7$ & 41 & 44 & 61.4 & 23.8 & 52.8 & 11.0 \\
    \midrule
    12 & No refresh & 0 & 1393 & 524 & 154 & 266 & 37.1 \\
     & Periodic, $T_{\mathrm{refresh}}=30$ & 46 & 13 & 19.0 & 4.20 & 1.45 & 1.11 \\
     & Periodic, $T_{\mathrm{refresh}}=50$ & 27 & 43 & 56.7 & 14.1 & 10.2 & 4.14 \\
     & Periodic, $T_{\mathrm{refresh}}=100$ & 13 & 93 & 54.9 & 17.5 & 7.86 & 4.19 \\
     & Bound, $\tau=0.5$ & 52 & 17 & 20.1 & 5.22 & 1.84 & 1.34 \\
     & Bound, $\tau=0.6$ & 37 & 7 & 9.95 & 2.24 & 1.28 & 0.895 \\
     & Bound, $\tau=0.7$ & 28 & 48 & 56.1 & 14.1 & 10.8 & 4.26 \\
    \bottomrule
  \end{tabular}
\end{table}

\subsection{Runtime Comparison}
\label{subsec:finance-runtime}

Finally, we examine runtime, measured in wall-clock seconds per
incremental update, i.e.\ per new trading day; the time of a refresh is
included in the update that triggers it. For $k=5$, the mean
update time ranges from $1.3 \cdot 10^{-4}$\,s for
\texttt{no\_refresh} to $2.2 \cdot 10^{-4}$\,s for
\texttt{refresh\_30}, i.e.\ one to two tenths of a
millisecond per day. About 30\% of the time without refresh is the copy
of the dense return matrix that each append makes
(Section~\ref{subsec:metrics}), an $\mathcal{O}(m_tn_t)$ cost that grows
with the number of days.

Figure~\ref{fig:finance-runtime} shows how this cost scales with the
working rank $k$ and the refresh cadence. In all cases, runtime grows
approximately linearly with $k$ (for \texttt{no\_refresh} from
$1.1\cdot 10^{-4}$\,s at $k=3$ to $2.7\cdot 10^{-4}$\,s at $k=20$), while
more frequent refreshes shift the curves upward by a nearly constant
amount (\texttt{refresh\_30} adds $0.8$--$0.9\cdot 10^{-4}$\,s at every $k$). The near-linear trend in
Figure~\ref{fig:finance-runtime} is consistent with the small-core
update structure discussed in
Sections~\ref{subsec:row-col-updates} and \ref{subsec:rank1-entry-updates},
together with the overall complexity discussion in
Section~\ref{subsec:metrics}. In practical terms, this means that
increasing the refresh frequency improves accuracy and stability, but
only at a moderate computational cost.

At this scale, a one-off full SVD of the final $1{,}643 \times 60$
returns matrix is also fast (a few milliseconds, about 30 times the
time of one incremental update with $k=5$ without refresh), but this
comparison is misleading. A single offline
decomposition solves only one static problem, whereas a live risk
system must absorb new observations continuously. Recomputing a full
SVD every time new data arrive would multiply that cost by the number
of updates and would quickly become prohibitive as either the asset
universe or the update frequency grows.

By contrast, the incremental SVD approach has a small and predictable
per-update cost, while periodic refresh offers a practical mechanism
for trading off speed against statistical fidelity. Combined with a
reasonable refresh schedule, the method therefore delivers near-oracle
covariance and portfolio-risk estimates while remaining suitable for
genuinely streaming environments, \emph{\textbf{scaling
naturally to high-frequency updates}}.

\section{Discussion}
\label{sec:discussion}

\subsection{Summary of Findings}

Across synthetic and real-data experiments, our results paint a consistent picture of when incremental SVD is accurate, stable, and computationally attractive.

On synthetic streams (Section~\ref{sec:synthetic}), the incremental updates behave exactly as the linear algebra suggests.
In the rank-1 streaming setting, \emph{no-refresh} runs accumulate reconstruction error and principal-angle drift over long horizons, while periodic refreshes (every $10^3$ and $5 \cdot 10^3$ updates) keep both the error ratio $\texttt{frob\_ratio}(t)$ from \eqref{eq:frob-ratio} and the angle to the optimal subspace $\texttt{angle\_to\_opt}$ defined in \eqref{eq:angle-metrics-opt} tightly bounded.
This is a direct consequence of the update structure: each rank-1 perturbation is projected back into a fixed rank-$k$ subspace, so any energy that repeatedly falls outside that subspace must eventually be ``forgotten'' unless we realign with a full SVD.
The a posteriori bound $\beta_t$ of Proposition~\ref{prop:error-control} holds at every step of these runs, and refreshing when $\beta_t/\|A_t\|_F$ exceeds $\tau$ (Section~\ref{subsec:bound-refresh}) is practically tied with periodic refresh at an equal number of full SVDs, while smaller $\tau$ trades more full SVDs for a lower error ratio.
The experiments also confirm the expected trade-off in the truncation rank $k$:
for $k \in \{5, 8, 12\}$ with true rank $5$, larger $k$ reduces the reconstruction error but increases the error ratio and the per-update cost (by a factor of 1.5 from $k=5$ to $k=12$ on that $50\times40$ stream), and the extra directions fit noise, so the angle to the optimal subspace grows towards $\tfrac{\pi}{2}$ (Section~\ref{subsec:k-effect}).
Row and column growth experiments show that Brand-style updates remain numerically stable as the matrix shape changes, provided $k$ is not too small relative to the intrinsic rank.
Scaling experiments show that, for fixed $k$, the rank-1 update time does not grow from the $50\times40$ to the $500\times400$ matrix, and that on the $500\times400$ matrix it grows faster than linearly with $k$ (a factor of 3.3 from $k=20$ to $k=40$); a full SVD of that matrix takes about 1400 times as long as one update with $k=5$.

In the financial application (Section~\ref{sec:finance}), we treat daily log returns of a 60-ETF universe as a time$\times$asset matrix and use incremental SVD as an \emph{online factor model} and \emph{low-rank covariance engine}.
Starting from an initial window $T_0$ of 250 days and then streaming about 1{,}400 additional days, we compare periodic refresh policies and truncation ranks $k \in \{3,5,8,12,20\}$, and run the bound-triggered rule for $k\in\{5,8,12\}$.
For moderate ranks ($k=5$--$12$) and refresh periods around 30 trading days, the incremental factors remain extremely close to the full-SVD factors: the final principal angle between subspaces is on the order of $10^{-2}$ radians, the relative covariance error is $4$--$6\cdot 10^{-4}$, and the relative ex-ante risk error is $1.5$--$4.7\cdot 10^{-6}$ for the equal-weight portfolio and $1.1$--$3.9\cdot 10^{-4}$ on average over random long--short portfolios (Table~\ref{tab:finance-bound-refresh}).
No-refresh runs exhibit large subspace drift (final angles 0.52--0.97 radians for $k=5$--$20$); at $k=3$ the no-refresh angle is small (0.084 radians) but the equal-weight risk error is the largest in the grid ($7.4\cdot 10^{-4}$). Very aggressive refreshing improves accuracy at the cost of more full-SVD calls.
The bound-triggered rule gives final-day errors comparable to periodic refresh; they follow the number of days since the last full SVD more closely than the policy, and because the floor $E_{\text{opt}}(t)/\|A_t\|_F$ is high on this stream, $\tau$ mainly sets the refresh budget.
A notable feature is the spike in risk approximation error around the 2020 market crash: an abrupt volatility regime shift creates large, correlated moves that neither full nor incremental low-rank models can forecast from pre-crash data; in our metrics this appears as a rise in covariance and risk errors in the weeks after the crash, which lasts until the next refresh (five days after the peak for refresh every 20 days, until late August 2020 for refresh every 100 days); without refresh the factor angle stays above 0.97 radians for the rest of the sample.

Overall, the experiments confirm that:
(i) incremental SVD can closely track the optimal rank-$k$ structure under moderate drift,
(ii) principal-angle and error-ratio metrics give interpretable diagnostics of when the approximation is degrading, and
(iii) with sensible choices of $k$ and refresh frequency, the method provides a practical accuracy--latency trade-off in 
financial settings.

\subsection{Refresh Policies and Adaptive Rank}

Section~\ref{sec:policies} introduced a design space of refresh and rank-selection policies; the experiments instantiate two refresh policies, periodic and bound-triggered, and support several qualitative conclusions.

First, \emph{periodic refresh} is empirically robust and easy to reason about.
In the synthetic rank-1 and mixed streams, refreshing every $10^3$ and every 200 updates, respectively, keeps $\texttt{frob\_ratio}$ close to $1$ (Sections~\ref{subsec:rank1} and~\ref{subsec:mixed}); on the rank-1 stream the angle to the optimal subspace stays small even without refresh (at most 0.044 radians over the seeds of Table~\ref{tab:bound-refresh}), and it is on the ETF stream without refresh that the factor angle drifts towards $\tfrac{\pi}{2}$ (Section~\ref{subsec:finance-factor}).
In the ETF experiments, refreshing every 20--100 trading days offers a clean bias--variance trade-off:
short refresh periods (e.g.\ 20--30 days) deliver almost indistinguishable factors and covariances compared to full SVD at every date, while less frequent refreshes (e.g.\ 100 days or never) gradually accumulate subspace error but save full-SVD calls.
These behaviors match the intuition that refreshes act as occasional ``re-projections'' onto the true dominant subspace, resetting the error that incremental updates cannot correct within a fixed-rank approximation.

Second, the bound-triggered rule \eqref{eq:bound-refresh-rule} replaces the fixed period by the a posteriori bound $\beta_t$ of Proposition~\ref{prop:error-control}; the test costs $\mathcal{O}(1)$ per step.
On the rank-1 stream (Section~\ref{subsec:bound-refresh}) the bound holds at every step up to roundoff but is loose without refresh ($\beta_t/E_{\text{inc}}(t)$ is about 80 after $10^4$ updates); refreshes reset it, and the median ratio drops to 7.6 with $T_{\mathrm{refresh}}=1000$ and to 1.2 with $\tau=0.05$.
At an equal number of full SVDs, $\tau=0.4$ and $T_{\mathrm{refresh}}=1000$ are practically tied, and smaller $\tau$ buys a lower error ratio with more full SVDs.
In the ETF experiments the floor $E_{\text{opt}}(t)/\|A_t\|_F$ is high (up to 0.45 for $k=5$), so $\tau$ acts as a refresh budget, and the final-day errors follow the number of days since the last full SVD more closely than the policy.
On neither stream does the rule beat periodic refresh at a matched number of full SVDs; what it adds is the guarantee $E_{\text{inc}}(t)\le\tau\|A_t\|_F$ after every step, as long as $\tau$ exceeds the floor $E_{\text{opt}}(t)/\|A_t\|_F$ (on one seed of the rank-1 stream the floor ends above $\tau=0.05$, and the guarantee is lost).
A guarantee on $\texttt{frob\_ratio}(t)$ \eqref{eq:frob-ratio} itself needs the trigger $\beta_t/E_{\text{opt}}(t_0)$ of Section~\ref{subsec:error-control}, which keeps $\texttt{frob\_ratio}(t)\le\gamma$ on append-only streams.

Third, the principal-angle metrics $\texttt{angle}$ (vs.~the initial basis) and $\texttt{angle\_to\_opt}$ (vs.~the instantaneous optimal subspace, when available in synthetic experiments) behave like \emph{concept-drift indicators}.
In both synthetic and financial streams we observe a tight coupling between large jumps in these angles and periods of structural change (e.g.\ accumulated noise, regime shifts).
This suggests angle-based refresh rules of the form ``refresh when $\texttt{angle\_to\_opt}$ exceeds a threshold'' or, in real systems where the optimal subspace is not available, ``refresh when the angle between successive bases or between new observations and their projections becomes too large.''

Fourth, the experiments illustrate how truncation rank $k$ interacts with refresh.
In the synthetic rank-effect experiment we choose $k$ at and above the true rank (5, 8, 12) to probe over-parameterization; the row and column streams (Section~\ref{subsec:structural}, $k=4$ against true rank 6) are under-parameterized.
In the financial setting we sweep $k\in\{3,5,8,12,20\}$ to span from very compact factor structures (3 factors, reminiscent of market/size/value-style models) to more saturated ones.
Higher rank increases the per-update time but does not systematically reduce the final errors against the rank-$k$ oracle: from $k=5$ upward the covariance error mostly falls, the factor angle does not, and $k=3$ gives the smallest angles (Table~\ref{tab:finance-bound-refresh} and Section~\ref{subsec:finance-risk}).
Taken together, these results motivate adaptive rank rules based on explained-variance thresholds (EVR) or residual ``novelty'' signals:
increase $k$ when the residual energy or new-row angle repeatedly exceed a tolerance, and decrease $k$ when EVR is far above the target and trailing singular values remain negligible.

We do not fully implement learned or fully adaptive policies in this work, but the metrics and behaviors we observe are precisely the signals a controller would need to \emph{automatically} decide when to refresh and how to adjust rank in a production system.

\subsection{Comparison with Alternative Approaches}

Our Brand-style incremental SVD is one of several paradigms for online low-rank modeling; each occupies a different point in the accuracy--efficiency--generality design space.
\emph{Randomized SVD}~\cite{halko2011randomized} is the method of choice when a full matrix snapshot is available and a fast batch recomputation suffices. In our setting it serves naturally as the ``refresh'' primitive; between refreshes, however, it offers no mechanism for absorbing individual updates.
\emph{Frequent Directions}~\cite{liberty2013simple,ghashami2016frequent} provides deterministic, space-optimal streaming sketches with provable spectral guarantees, but produces a sketch matrix rather than an explicit $(U,\Sigma,V)$ decomposition, making it less suited for applications that require singular vectors (e.g.\ factor loadings in risk models) or entry-level updates.
\emph{Stochastic subspace trackers} such as GROUSE~\cite{balzano2010grouse} and Oja's rule~\cite{oja1982simplified} are highly scalable and handle incomplete observations, but converge only asymptotically and do not yield singular values, limiting their use for covariance estimation or explained-variance diagnostics; see~\cite{balzano2018streaming} for a unified survey of these subspace-tracking paradigms.
For the financial covariance application specifically, the main alternatives are parametric models such as DCC-GARCH~\cite{engle2002dcc}, which directly model time-varying correlations but scale quadratically in the number of assets, and shrinkage estimators~\cite{ledoit2004wellconditioned} or thresholded factor models (POET)~\cite{fan2013large}, which produce well-conditioned covariance matrices but are typically applied in a rolling-window batch fashion rather than updated incrementally.
Our experiments show that incremental SVD with periodic or bound-triggered refresh achieves near-oracle covariance and risk accuracy while remaining fully incremental; a systematic empirical comparison with these alternatives on larger asset universes is an interesting direction for future work.

\subsection{Limitations and Scope}

Our study has several limitations that we view as opportunities for future work. 

\paragraph{Scale and sparsity.}
All experiments use a dense implementation of SVD and incremental updates.
This is deliberate, as it allows us to measure errors exactly and to run full-SVD baselines at each step in synthetic settings, but it restricts us to moderate matrix sizes.
Real recommender systems and factor models operate on sparse and extremely large matrices; extending our framework to sparse, GPU-accelerated, or randomized SVD variants is an important next step. 

Recent work shows that truncated-SVD updates can be implemented efficiently in the sparse regime while preserving representation quality~\cite{deng2024fasttsvd}.
Complementarily, projection-based updating schemes that reduce the dimension of the update subspace suggest a practical route to handling larger update blocks more efficiently~\cite{vecharynski2014fastlsi}.

\paragraph{Data modeling choices.}
In the financial application, we work with daily log returns of a 60-ETF universe over a few years.
This is realistic enough to exhibit volatility regimes (e.g.\ the 2020 crash) but far from high-frequency order-book data or cross-sectional universes with thousands of stocks; we also ignore transaction costs, constraints, and turnover, focusing solely on ex-ante risk.

\paragraph{Policy space coverage.}
We instantiate two refresh policies with fixed ranks: periodic refresh, with a grid of $(k, \texttt{refresh\_every})$ combinations in the ETF experiments, and the bound-triggered rule \eqref{eq:bound-refresh-rule}, on the rank-1 stream and for $k\in\{5,8,12\}$ on the ETF stream.
The bound-triggered rule controls the error relative to $\|A_t\|_F$; the trigger $\beta_t/E_{\text{opt}}(t_0)$ of Section~\ref{subsec:error-control}, which bounds $\texttt{frob\_ratio}$ on append-only streams, is not run in our experiments.
The ETF comparison uses final-day snapshots, which depend strongly on the days since the last full SVD; a time-averaged comparison, as in Section~\ref{subsec:bound-refresh}, would separate the policies more fairly.
Angle-based refresh and truly adaptive rank selection are formalized in Section~\ref{sec:policies} and supported by our metrics, and we leave their evaluation to future work.
Likewise, we do not attempt to learn refresh or rank policies from data (e.g.\ via reinforcement learning or bandit-style adaptation).

\paragraph{Theoretical analysis.}
Our work is primarily empirical and algorithmic. Proposition~\ref{prop:error-control} gives a computable a posteriori bound on the error accumulated between refreshes, but we do not provide a priori bounds on error accumulation, nor bounds on subspace drift, under streaming updates.
The a posteriori bound is also loose on long segments without refresh (about $80\times$ the true error after $10^4$ rank-1 updates), because the triangle inequality behind it ignores cancellation between the discarded parts.
The observed coupling between principal angles, error ratios, and refresh decisions suggests that such bounds may be achievable for certain stochastic models, but we leave this for future theoretical work.
Recent work in a different domain (continual learning of classifiers from pre-trained features) proves via a recurrence relation that continually truncating an appropriate fraction of SVD factors keeps error bounded over long task sequences~\cite{peng2025loranpac}, suggesting that similar truncation-based bounds may be achievable for streaming row/column updates as well.

Despite these limitations, the experiments demonstrate that a relatively simple incremental SVD engine, combined with a metric-rich monitoring layer, already provides actionable insights and performance guarantees that are relevant for dynamic recommender systems and online factor/risk models.

\section{Conclusion}
\label{sec:conclusion}

We studied incremental SVD as a practical tool for maintaining low-rank approximations of dynamic matrices that evolve through row additions, column additions, and rank-1 entry updates.
Starting from a Brand-style algebraic core, we built an incremental engine that supports bounded-rank updates, optional refreshes via full SVD, and a suite of diagnostics that quantify reconstruction error, subspace drift, and computational cost over time.

On the methodological side, the paper contributes:
(i) a unified incremental SVD framework covering heterogeneous update types;
(ii) a computable a posteriori bound on the error accumulated between refreshes (Proposition~\ref{prop:error-control});
(iii) a metric-based evaluation layer combining error ratios, principal angles, explained-variance ratios, and runtime; and
(iv) a refresh and rank-selection design space including periodic, error-based, and angle-based policies, as well as EVR- and novelty-driven rank adaptation.
On the empirical side, synthetic streaming experiments isolate the roles of rank, drift, and refresh frequency and compare periodic refresh with a bound-triggered rule built on Proposition~\ref{prop:error-control}; 
a financial case study demonstrates that a low-rank incremental factor model can track full-SVD covariances and portfolio risks in an online, time-series setting.

A recurring theme across these domains is that \emph{monitoring the right quantities} is as important as the update formulas themselves.
Principal angles provide an interpretable notion of subspace drift and concept shift; error ratios benchmark the incremental approximation against the best possible rank-$k$ model; and simple timing metrics compare the cost of an incremental update with that of a full recomputation.
Our experiments indicate that, with sensible choices of rank and refresh frequency, incremental SVD offers a useful accuracy--latency trade-off in 
financial applications and that it \emph{scales naturally to high-frequency updates} in scenarios where repeated full SVDs are infeasible.

There are several promising directions for future work.
On the algorithmic side, extending the framework to sparse settings, and combining it with randomized or streaming SVD variants, would make the method directly applicable to industrial-scale 
enterprise risk engines.
On the policy side, learning refresh and rank strategies from data---for example via reinforcement learning or online convex optimization---could turn our metric signals into fully automated controllers.
Finally, applying the same ideas to other dynamic-matrix domains (e.g.\ sensor networks, robotics, digital-twin environments, or high-frequency limit-order books) would further test the robustness and usefulness of incremental SVD as a general-purpose building block for online low-rank modeling.

\begin{appendices}

\section{Detailed Incremental Update Routines}
\label{app:updates}

This appendix expands Algorithm~\ref{alg:incrementalsvd} by giving detailed pseudocode for each update primitive
used by the \texttt{IncrementalSVD} engine.

\setcounter{algorithm}{0}
\renewcommand{\thealgorithm}{A\arabic{algorithm}}
\renewcommand{\theHalgorithm}{A\arabic{algorithm}}% distinct hyperref anchors from Algorithm 1

\begin{algorithm}[H]
\caption{\texttt{RowAppendUpdate} (detailed)}
\label{alg:rowappend}
{\footnotesize
\begin{algorithmic}[1]
\Require $U,\Sigma,V^\top$ (rank $r\le k$), new row $x^\top$, tolerance \texttt{tol}
\State $p \gets x^\top V$
\State $z \gets x^\top - pV^\top$, \quad $\rho \gets \lVert z\rVert_2$
\If{$\rho > \texttt{tol}$}
  \State $K \gets \mbox{\scriptsize $\left[\!\! \begin{array}{rr} \Sigma & 0 \\ p & \rho \end{array} \!\!\right]$}$, \quad $\bar V^\top \gets \mbox{\scriptsize $\left[\!\! \begin{array}{c} V^\top \\ z/\rho \end{array} \!\!\right]$}$
\Else
  \State $K \gets \mbox{\scriptsize $\left[\!\! \begin{array}{c} \Sigma \\ p \end{array} \!\!\right]$}$, \quad $\bar V^\top \gets V^\top$ \Comment{drop $z$, no new direction}
\EndIf
\State $(U_K,\Sigma',V_K^\top) \gets \texttt{SVD}(K)$ \Comment{thin SVD}
\State $U' \gets \mbox{\scriptsize $\left[\!\! \begin{array}{rr} U&0\\0&1 \end{array} \!\!\right]$}U_K$
\State $V'^\top \gets V_K^\top \bar V^\top$
\State \Return $(U',\Sigma',V'^\top)$ (then truncate to $k$ if needed)
\end{algorithmic}
}
\end{algorithm}

\begin{algorithm}[H]
\caption{\texttt{ColAppendUpdateViaTranspose} (detailed)}
\label{alg:colappend}
{\footnotesize
\begin{algorithmic}[1]
\Require $U,\Sigma,V^\top$, new column $y$, tolerance \texttt{tol}
\State Treat as row-append on $A^\top$: $\widetilde{U} \gets V$, \, $\widetilde{\Sigma}\gets\Sigma$, \, $\widetilde{V}^\top \gets U^\top$
\State $(\widetilde{U}',\Sigma',\widetilde{V}'^\top) \gets \texttt{RowAppendUpdate}(\widetilde{U},\widetilde{\Sigma},\widetilde{V}^\top,y^\top;\texttt{tol})$
\State Transpose back: $U' \gets (\widetilde{V}'^\top)^\top$, \, $V'^\top \gets (\widetilde{U}')^\top$
\State \Return $(U',\Sigma',V'^\top)$ \, (then truncate to $k$ if needed)
\end{algorithmic}
}
\end{algorithm}

\begin{algorithm}[H]
\caption{\texttt{RankOneUpdate} (detailed)}
\label{alg:rank1}
{\footnotesize
\begin{algorithmic}[1]
\Require $U,\Sigma,V^\top$ (rank $r\le k$), index $(i,j)$, increment $\delta$
\State $s_u \gets (U_{i,:})^\top$, \, $s_v \gets (V^\top)_{:,j}$
\State $K \gets \Sigma + \delta \, s_u s_v^\top$
\State $(U_K,\Sigma',V_K^\top) \gets \texttt{SVD}(K)$
\State $U' \gets U U_K$, \, $V'^\top \gets V_K^\top V^\top$
\State \Return $(U',\Sigma',V'^\top)$ \, (rank preserved; truncate if needed)
\end{algorithmic}
}
\end{algorithm}

\end{appendices}

\backmatter

\section*{Statements and Declarations}

\subsection*{Data availability}
Daily prices of 60 US-listed exchange-traded funds were obtained from Yahoo Finance through the open-source Python package \texttt{yfinance}.
The tickers are listed in the script \texttt{experiments/prepare\_etf\_data.py} in the code repository (see Code availability).
The script downloads adjusted closing prices (adjusted for splits and dividends) for 1~January 2010 to 31~December 2024 and keeps only the trading days on which all 60 prices are available; this leaves 1644 trading days, from 19~June 2018 to 30~December 2024.
The experiments convert these prices to daily log-returns and center each asset as described in Section~\ref{subsec:finance-data}.
The data are provided by a third party, so the author does not redistribute them; the dataset can be rebuilt with \texttt{experiments/prepare\_etf\_data.py}.
The synthetic data are generated by the released code with fixed random seeds.

\subsection*{Code availability}
The code is available at \url{https://github.com/stilyanst/iSVD}.
It contains the Python package \texttt{isvd} (incremental SVD, refresh policies and metrics; it depends only on NumPy) and the experiment scripts: \texttt{experiments/reproduce.py} regenerates all figures and tables of this paper, and \texttt{requirements.txt} lists the exact package versions used.
The full-SVD timings in Sections~\ref{subsec:scaling} and~\ref{subsec:finance-runtime} and the share of the matrix copy in the update time (Section~\ref{subsec:finance-runtime}) were measured separately, in the timing setup of Section~\ref{subsec:metrics}, and are not produced by \texttt{experiments/reproduce.py}.
The code is released under the MIT license.
Version 1.0.0, used for this paper, is archived on Zenodo~\cite{staykov2026isvd}, \url{https://doi.org/10.5281/zenodo.24481009}.

\bibliography{reference}

\end{document}